\documentclass[12pt]{article}
\usepackage{amsmath,amsthm,amssymb,amsfonts}

\newfont{\lie}{eufm10 at 12pt}
\newfont{\liepequenos}{eufm10 at 10pt}

\newfont{\corpos}{msbm10 at 12pt}
\newfont{\corpospequenos}{msbm10 at 10pt}

\newcommand{\Cpequeno}{\mbox{\corpospequenos \symbol{67}}}   
\newcommand{\Ppequeno}{\mbox{\corpospequenos \symbol{80}}}    
\newcommand{\Rpequeno}{\mbox{\corpospequenos \symbol{82}}}    
\newcommand{\Tpequeno}{\mbox{\corpospequenos \symbol{84}}}    

\newcommand{\C}{\mbox{\corpos \symbol{67}}}          
\newcommand{\Pro}{\mbox{\corpos \symbol{80}}}        
\newcommand{\R}{\mbox{\corpos \symbol{82}}}          
\newcommand{\T}{\mbox{\corpos \symbol{84}}}          

\newcommand{\g}{\mbox{\lie g}}       %
\newcommand{\h}{\mbox{\lie h}}       %
\newcommand{\gl}{\mbox{\lie gl}}     
\newcommand{\m}{\mbox{\lie m}}       %

\newcommand{\uni}{\mbox{\lie u}}     %

\newcommand{\symp}{\mbox{\lie sp}}

\newcommand{\DCON}{\mbox{\lie D}}   
\newcommand{\XIS}{\mbox{\lie X}}   

\newcommand{\co}{connection}

\newcommand{\rdoisn}{\R^{2n}}                    %
\newcommand{\jm}{{\cal J}(M)}           
\newcommand{\jnab}{{\cal J}^{\nabla}}             
\newcommand{\hnab}{{\cal H}^{\nabla}}             
\newcommand{\zl}{{\cal Z}^l}           
\newcommand{\zo}{{\cal Z}^0}           

\newcommand{\nab}[2]{{\nabla_{_{#1}}#2}}
\newcommand{\Tr}[1]{{\mathrm{Tr}}\,#1}

\newcommand{\End}[1]{{\mathrm{End}}\,#1}

\newcommand{\Id}{{\mathrm{Id}}}
\newcommand{\dx}{{\mathrm{d}}}
\newcommand{\db}{\overline{\partial}}
\newcommand{\cinf}[1]{{\mathrm{C}}^\infty_{#1}}

\newtheorem{teo}{Theorem}[section]
\newtheorem{lema}{Lemma}[section]
\newtheorem{coro}{Corollary}[section]
\newtheorem{prop}{Proposition}[section]

\numberwithin{equation}{section}

\begin{document}

\thispagestyle{empty}
\begin{center}
\vspace{2cm} {\huge{\bf{Twistor Theory of Symplectic \vspace{.7cm}

Manifolds}}}

\vspace{2cm} {\sc{R. Albuquerque}}

{\small{rpa@uevora.pt}}

Departamento de Matem\'atica

Universidade de \'Evora

7000 \'Evora

Portugal

\vspace{.9cm} {\sc{J. Rawnsley}}

{\small{j.rawnsley@warwick.ac.uk}}

Mathematics Institute

University of Warwick

Coventry CV47AL

England

\vspace{1.2cm}

November 2004

\small{ Abstract

This article is a contribuition to the understanding of the
geometry of the twistor space of a symplectic manifold. We
consider the bundle ${\cal Z}^l$ with fibre the Siegel domain
$Sp(2n,R)/U(n)$ existing over any given symplectic manifold $M$.
Then, while recalling the construction of the celebrated almost
complex structure induced on ${\cal Z}^l$ by a symplectic
connection on $M$, we study and find some specific properties of
both. We show a few examples of twistor spaces, develop the
interplay with the symplectomorphisms of $M$, find some results
about a natural almost-hermitian structure on ${\cal Z}^l$ and
finally discuss the holomorphic completeness of the respective
``Penrose transform".}

\end{center}
\vspace{2.5cm}

Let $(M,\omega)$ be a smooth symplectic manifold of dimension
$2n$. Then we may consider the bundle
\begin{equation}
\pi:{\cal Z}^l\longrightarrow M              \label{eq0}
\end{equation}
of all complex structures $j$ on the tangent spaces to $M$
compatible with $\omega$. Having fibre a cell, the bundle becomes
interesting if it is seen with a particular and well known almost
complex structure, denoted $\jnab$, --- with which we start to
treat ${\cal Z}^l$ by the name of ``Twistor Space'' of the
symplectic manifold $M$. The almost complex structure is induced
by a symplectic connection on the base manifold and its
integrability equation has already been studied
(\cite{Nannicini1,Nannicini2,Rawnsley3,Vais2}).

The purpose of this article is two-fold. First, we wish to present
results on the complex geometric nature of the twistor space and,
second, to show some of the applications to the study of
symplectic connections. The almost complex structure $\jnab$ is
defined in a very peculiar way, in the sense, for example, that no
matter which complex structure is given on the base space $M$ the
bundle projection is never pseudo-holomorphic. The complex
geometry of the twistor space seems, up to a certain degree which
we compute, much more difficult to understand than that one of a
regular holomorphic fibre bundle. This is particularly true in the
study of the Penrose Transform, which we see as a direct image
from cohomology of complex analytic sheaves on $\zo$ to real
$\cinf{}$ sheaves on $M$. We obtain a ``vanishing theorem'' which
is what one would expect if both spaces were complex and $\pi$
holomorphic (section 6). We will also show examples of twistor
spaces and a compact generalized ``omega'' twistor space, in
section 3, and give results about a natural Riemannian structure
on ${\cal Z}^0$ in section 5 (following a different approach, such
structure has been considered in \cite{Nannicini1}).

Concerning the second purpose of this work, we relate the twistor
almost complex structure to the action of the group of
symplectomorphisms of $M$ on the space of symplectic connections.
In passing by this independent subject, we present our methods in
its study --- and apply them on the case of translation invariant
symplectic connections, finally just to find a known result proved
by other means, {\it cf.} \cite{Rawnsley2}. Later we give a new
criteria to decide when are two given symplectic connections the
affine transformation of one another via a given symplectomorphism
of $M$. We may then claim to have found a description of all germs
of twistor spaces of a Riemann surface (sections 2,3,4). Moreover,
since the integrability equation is immediately and always
satisfied in case $n=1$, with the Levi-Civita connection in
particular, we describe a new $\cinf{}$ sheaf canonically defined
on $M$, arising from a $\C$-analytic sheaf on $\zo$. With this we
hope to have contributed to future studies in the field of
twistors.

\vspace{1cm}

\section{The structure of twistor space}

Throughout the text let $G=Sp(2n,\R)$ and $U^l=U(n-l,l)$, {\it
ie.} respectively the groups of symplectic and pseudo-unitary
transformations of $\R^{2n}$. Also let $\g$ and $\uni^l$ denote
their respective Lie algebras.

Let $\omega=\sum_{i=1}^n\dx x^i\wedge\dx y^i$. Recall that the
complex symmetric space $G/U^l$ consists of all {\it real}
endomorphisms $J$ such that
\begin{eqnarray*}
\begin{cases}
 J^2=-\Id,\ \ \ \ \ \omega=\omega^{1,1}\ \mbox{for $J$},
\\
\omega(\ ,J\ )\ \mbox{has signature}\ (2n-2l,2l) .
\end{cases}
\end{eqnarray*}
It is known that these spaces are biholomorphic to open cells of
the flag manifold $Sp(n)/U(n)$ and that the latter lyes on the
complex grassmannian of $n$-planes in $\C^{2n}$ as the space of
lagrangians. Moreover, the boundaries of the $G/U^l$ do not give
rise to complex structures under the inverse of the map $J\mapsto$
$\sqrt{-1}$-eigenspace of $J$, {\it cf.} \cite{Tese}. In a word,
there are no other complex structures $J$ of $\R^{2n}$ for which
$\omega$ is type $(1,1)$.

The vector space $T_J\,G/U^l$ identifies with
\[ \m_J=\left\{A\in\g:\ AJ=-JA \right\} =[\g,J] \]
so that $\g=\m_J+\uni^l$ and the complex structure is given by
left multiplication by $J$. We have
\begin{equation}
  [\g,\m_J]=\m_J,\ \ \ \ \ [\m_J,\m_J]\subset\uni^l,\ \ \ \
[\uni^l,\uni^l]\subset \uni^l        \label{eq2}
\end{equation}
as one can easily compute. It is known that all the symmetric
spaces $G/U^l$ are biholomorphic to the first one, when $l=0$,
which is the same as the Siegel Domain or Siegel Upper Half Space
$\{X+iY\in\C^{\frac{n(n+1)}{2}}: \ X,Y \mbox{ symmetric matrices,
$Y$ positive definite}\}$, and hence that they are all Stein
manifolds. Up to this moment we have met the $n+1$ possible
connected fibres, according to $0\leq l\leq n$, which will
correspond to the various twistor spaces of a symplectic manifold.

Now consider the manifold $(M,\omega)$ and the twistor bundle
\eqref{eq0} with fibre $G/U^l$. When necessary we denote this
space by ${\cal Z}^l_M$. We have a short exact sequence
\begin{equation}
0\longrightarrow{\cal V}\longrightarrow T{\cal
Z}^l\stackrel{\dx\pi}\longrightarrow \pi^*TM\longrightarrow0
\label{eq3}
\end{equation}
of vector bundles over the manifold ${\cal Z}^l$. If we let
$E=\pi^*TM$ and denote by $\Phi$ the canonical section of
$\End{E}$, given by
\[   \Phi_j=j ,\]
then we may identify ${\cal V}=\ker\,\dx\pi$ with
$[\symp(E,\omega),\Phi]$. This is justified by the theory above
applied to the fibre $\pi^{-1}(x)$, which is thus a complex
manifold, for all $x\in M$.

We use now a symplectic connection $\nabla$ on $TM$ (see section
2), in order to construct a horizontal distribution ${\cal
H}^\nabla$ which splits the sequence \eqref{eq3}. The following
result from the theory of twistor spaces is adapted to our
situation, so we recall the proof briefly --- to introduce
notation and for later reference.
\begin{prop}[\cite{Obri}] \label{prop1.1}
\ ${\cal H}^\nabla =\{ X\in T{\cal Z}^l:\
(\pi^*\nabla)_{_X}\Phi=0 \}$ is a complement for $\cal V$ in
$T{\cal Z}^l$.
\end{prop}
Let $F^s(M)$ be the symplectic frames bundle of $M$, consisting of
all linear symplectic isomorphisms
\[  p:\rdoisn\longrightarrow T_xM.      \]
This is a principal $G$ bundle over $M$ and $\nabla$ is determined
by the $\g$-valued 1-form $\alpha$ on $F^s(M)$, given by
\[   \nab{X}{(sv)}=s(s^*\alpha)(X)v         \]
for any section $s$ of $F^s(M)$, for $X\in \Gamma TM,\ v\in
\rdoisn$, and
 by $\alpha(\tilde{A})=A$
where $A\in\g$ and $\tilde{A}$ is the vector field
\[    \tilde{A}_p={\frac{\dx}{\dx t}}_{|_{0}}\ p\circ \exp tA.        \]
$\ker\alpha$ is a horizontal distribution on $F^s(M)$. Fixing
$J_0\in G/U^l$, we get a map
\begin{equation}
\pi_1:F^s(M)\longrightarrow {\cal Z}^l
\label{pium}
\end{equation}
given by
\[         \pi_1(p)=pJ_0p^{-1}.      \]
The derivative of $\pi_1$ maps the horizontal distribution on
$F^s(M)$ onto a ho\-rizontal distribution on ${\cal Z}^l$. The
proposition is proved by showing that this distribution coincides
with ${\cal H}^\nabla$.

Notice \eqref{pium} is a principal $U^l$ bundle over ${\cal Z}^l$
and that $E$ is also associated to $\pi_1$, due to reduction.
Since $\Phi$ is a section of $\End{E}$ corresponding with the
constant \textit{equivariant} function
$\hat{\Phi}(p)=p^{-1}\Phi_{\pi_1(p)}p=J_0$ on $F^s(M)$ and since
$(\pi^*\nabla)\Phi$ corresponds to the 1-form
\[   \dx J_0+[\alpha,J_0],     \]
one may then complete the proof of the proposition.

Denoting by
\[ P:T\zl\longrightarrow {\cal V} \]
the projection onto $\cal V$ with kernel $\hnab$, we have the
following consequences.
\begin{prop}[\cite{Obri}] \label{prop1.2}
(i) \,$(\pi^*\nabla)\Phi=[P,\Phi]=2P\Phi$.
\\
(ii) $T\zl=F^s(M)\times_{U^l}\m_{J_0}\oplus\rdoisn$.
\end{prop}

We may now define, preserving the direct sum $\hnab\oplus{\cal
V}$, the twistor almost complex structure
\[  \jnab=(J^h,J^v) ,\]
on each point $j\in\zl$ as follows: since $\dx\pi:\hnab\rightarrow
E$ is an isomorphism, we transport $j$ from $E$ to the horizontal
bundle. $J^h$ is thus, essentially, $j$ itself. The vertical part
$J^v$ consists of left multiplication with $j$, just like in the
Siegel Domain.

Let $i=\sqrt{-1}$ and let $j^+=\frac{1}{2}(1-ij)$,
$j^-=\frac{1}{2}(1+ij)$ be the projections from $TM\otimes\C$
onto, respectively, $T^+M$ and $T^-M$, for any $j\in\zl$. The
integrability equation for $\jnab$ follows from the next theorem,
a result which we present in greater generality for later
convenience. Notice $\jm$ is the general twistor space consisting
of all complex structures, {\it ie.} the bundle with fibre
$GL(2n,\R)/GL(n,\C)$, and that, of course, precisely in the same
lines of the case we have been considering, any linear connection
on $M$ defines a twistor almost complex structure on $\jm$.
\begin{teo}[\cite{Obri}]  \label{teo1.1}
Let $Z$ be an almost complex manifold and $\pi:Z\rightarrow M$ be
a smooth submersion onto $M$ with fibres which are smoothly
varying complex manifolds. Suppose that $Z$ has a horizontal
distribution ${\cal H}^Z$ which is $j$-related to the horizontal
distribution ${\cal H}^\nabla$ of a connection $\nabla$ on $TM$
via a pseudo-holomorphic smooth {\rm fibre preserving} map
\[     j:Z\longrightarrow \jm .      \]
Then integrability of $J^Z$ implies that the torsion $T$ and
curvature $R$ of $\nabla$ satisfy
\begin{equation}
     j^+T_x(j^-X,j^-Y)=0,\hspace{1.3cm}j^+R_x(j^-X,j^-Y)j^-=0  \label{eq4}
\end{equation}
for all $j\in j(Z)$ and $X,Y\in T_xM$. If $j$ is an immersion
these conditions are also sufficient.
\end{teo}
In the case of $\zl$, $j$ is just the inclusion map and equations
\eqref{eq4} are equivalent to the vanishing of the Weyl part of
the curvature of the given symplectic connection (see section 2,
\cite{Rawnsley3} and \cite{Tese} for definition and references).
In this regard, this Weyl tensor plays a role in symplectic
geometry identical to that Weyl tensor of a metric connection in
the theory of twistor spaces in Riemannian geometry, {\it cf.}
\cite{Atiyah,Obri}.

We remark that a proof of the independency of the equations
\eqref{eq4} from the signature, which varies with $l$, is also
written down in \cite{Tese}. Hence, either all the $\zl$ --- or
none --- have integrable almost complex structures.

The following is particular to the symplectic framework.
\begin{teo}  \label{teo1.2}
If ${\cal J}^{\nabla^1}={\cal J}^{\nabla^2}$ then
$\nabla^1=\nabla^2$.
\end{teo}
\begin{proof} Let $A=\nabla^2-\nabla^1$ and define $\underline{A}\in\Gamma S^3T^*M$ by
\begin{eqnarray*}
  \underline{A}(X,Y,Z) & = & \omega(A(X)Y,Z)\\
  & = & \omega(A(Y)X,Z)\ =\ \omega(A(Z)Y,X)
\end{eqnarray*}
for all $X,Y,Z\in TM$.

Now suppose that ${\cal J}^{\nabla^1}={\cal J}^{\nabla^2}$ and $u$
is any $\nabla^1$-horizontal vector field of type (1,0). Then
$u=u_2+v$, with $v$ vertical, also a (1,0) vector field, and $u_2$
a horizontal vector field for $\nabla^2$. Hence
\begin{eqnarray*}
 [v,\Phi]& =& \pi^*\nabla_u^2\Phi \\
 &=&[\pi^*A(\dx\pi\,u),\Phi]\ \in\, {\cal V}^{(1,0)}
\end{eqnarray*}
(by propositions \ref{prop1.1}, \ref{prop1.2} and noticing that
$[J_0,\m^+_{J_0}]\subset\m^+_{J_0}$). In the base manifold this
translates as
\[ [A_x(j^+X),j]\ \ \mbox{is}\ \ (1,0),\ \ \ \ \
\ \ \forall j\in\pi^{-1}(x),\ \forall X\in T_xM .\] Equivalently
this means the projection to ${\cal V}^c$ is (1,0), or
\[j^-A(j^+X)j^+=0,\ \ \ \ \ \ \ \forall j,\ \forall X . \]
The equality $\underline{A}^{3,0}=0,\ \forall j$, follows
immediately. This says $\underline{A}$ must take values in the
largest $G$-invariant subspace of symmetric tensors which satisfy
\begin{eqnarray*}
\underline{A}(J_0^+X,J_0^+Y,\,\ldots)=0
\end{eqnarray*}
for some fixed $J_0$ and all $X,Y,\,...$\,. Indeed, since any
$j=gJ_0g^{-1}$ for some $g\in G$, we will also have
$j^+=\frac{1}{2}(1-ij)=gJ_0^+g^{-1}$ and therefore
\begin{eqnarray*}
 0=(g^{-1}\cdot \underline{A})(J_0^+X,J_0^+Y,\,\ldots) & = &
          \underline{A}(gJ_0^+X,gJ_0^+Y,\,\ldots)\\
          & = &  \underline{A} (j^+gX,j^+gY,\,\ldots).
\end{eqnarray*}
But $S^k(\C^{2n})$ is irreducible under $G$ for all $k$, so
$\underline{A}$ must be $0$.
\end{proof}

\vspace{1cm}

\section{Symplectic connections}

\def\cyclic{\mathop{\kern0.9ex{{+}
\kern-2.2ex\raise-.28ex\hbox{\Large\hbox{$\circlearrowright$}}}}\limits}

Let $M,N$ be two manifolds and $\sigma:M\rightarrow N$ a
diffeomorphism between them. Let $\nabla$ be a linear connection
on $M$. Recall that we can define another connection on $N$ by
\[ \left(\sigma\cdot\nabla\right)_{_X}Y =\sigma\cdot\left(\nab{\sigma^{-1}\cdot
X}{\sigma^{-1}\cdot Y} \right),     \] where $X,Y\in\XIS_{N}$ and
where
\[     \sigma\cdot Z\,_y=\dx\sigma(Z_{\sigma^{-1}(y)})        \]
for any $Z\in\XIS_{M},\ y\in N$. It is well defined, at least, on
paracompact manifolds.

Indeed, from any tensor on $M$ we can define another one on $N$.
Notice as well that $\sigma\cdot fZ=(f\circ\sigma^{-1})\sigma\cdot
Z=\sigma\cdot f\ \sigma\cdot Z$, for all $f\in\cinf{M}$, so we
prove the last statement and check Leibniz rule for
$\sigma\cdot\nabla$. Furthermore, the torsion and curvature
satisfy
\begin{equation}
    T^{\sigma\cdot\nabla}=\sigma\cdot T^\nabla,\ \ \ \ \ \ \
    R^{\sigma\cdot\nabla}=\sigma\cdot R^\nabla ,                \label{eq6}
\end{equation}
since $\sigma\cdot [Z,W]=[\sigma\cdot Z,\sigma\cdot W]$. Obvious
composition rules are satisfied and
\begin{equation}
 {\left(\sigma\cdot\nabla\right)}_{_X}\omega=\sigma\cdot\left(\nab{\sigma^{-1}\cdot
 X}{\sigma^*\omega}\right)                         \label{eq7}
\end{equation}
for any form $\omega$ on $N$. For instance, let us prove the last
formula:
\begin{eqnarray*}
 \sigma\cdot\left(\nab{\sigma^{-1}\cdot X}{\sigma^*\omega}\right)(Y_1,\ldots,Y_q) &  =  &
 \left(\nab{\sigma^{-1}\cdot X}{\sigma^*\omega}\right)
 (\sigma^{-1}\cdot Y_1,\ldots,\sigma^{-1}\cdot Y_q)\ \ \ \ \ \ \
\end{eqnarray*}
\vspace{-.9cm}
\begin{eqnarray*}
  & \ \ \ \ \ \ =  & (\sigma^{-1}\cdot X)
  \left(\sigma^*\omega(\sigma^{-1}\cdot Y_1,\ldots,\sigma^{-1}\cdot Y_q)\right)  \\
  & &  \hspace{1.7cm} -\sum_i\sigma^*\omega\left(\sigma^{-1}\cdot Y_1,\ldots,
  \nab{\sigma^{-1}\cdot X}{\sigma^{-1}\cdot Y_i},\ldots,\sigma^{-1}\cdot Y_q\right) \\
    & \ \ \ \ \ \ = &   \dx(\omega_\sigma(Y_1,\ldots,Y_q))(\dx\sigma^{-1}(X))-
    \sum\omega\left(Y_1,\ldots,\left(\sigma\cdot\nabla\right)_{_X}Y_i,\ldots,Y_q\right)  \\
   & \ \ \ \ \ \ =  &   \left(\sigma\cdot\nabla\right)_{_X}\omega\, (Y_1,\ldots,Y_q).
\end{eqnarray*}
As we said before, $\sigma^{-1}\cdot\omega=\sigma^*\omega$.
\vspace{4mm}\\
{\bf Remark.} In a marginal outlook, the above may be applied in
Chern-Weil theory to find that all characteristic classes on $TM$,
induced from multilinear forms $f^i:\otimes^i\g\rightarrow\R$ or
$\C$, say $H$-invariant where $H$ is some Lie group, are fixed
points of cohomology for every diffeomorphism preserving some
$H$-structure of $M$. The proof goes as follows. Taking any
$H$-connection $\nabla$, assumed to exist, we then have
\begin{eqnarray*}
f^i(R^{\sigma^{-1}\cdot\nabla},\ldots, R^{\sigma^{-1}\cdot\nabla})
& = & f^i(\dx\sigma^{-1}(\sigma^*R^\nabla)\dx\sigma,
\,\ldots\,,\dx\sigma^{-1}(\sigma^*R^\nabla)\dx\sigma )\\
 & = & \sigma^*f^i(R^\nabla,\ldots, R^\nabla).
\end{eqnarray*}
Hence, by the independence of the induced de Rham cohomology
classes from the choice of the connection, the former are fixed
points for $\sigma^*$. Of course the result is
 interesting when $Dif\!f(M)$ has many arcwise-connected components.
\vspace{4mm}

From now on we are interested in the case where $M$ and $N$ are
symplectic manifolds and $\sigma$ is a symplectomorphism. Recall
that a linear connection on $(M,\omega)$ is called symplectic if
$\nabla\omega=0$ and if it is torsion free. In such case, by
formulae \eqref{eq6},\eqref{eq7}, we have that $\sigma\cdot\nabla$
is symplectic too. In particular we have an action
\[   Symp(M,\omega)\times {\cal A}\longrightarrow{\cal A}   \]
on the space of symplectic connections, which preserves the
subspace of flat connections. $\cal A$  is never empty.
\begin{teo}[P.Tondeur, {\it cf.} \cite{Cahen}]  \label{teo2.1}
Every symplectic manifold admits a symplectic connection.
\end{teo}
Furthermore, if a Lie group $H$ acts on $M$ by symplectomorphisms
and thus on the space of connections, then $M$ has a $H$-invariant
\co\ if and only if it has a $H$-invariant symplectic \co.

Notice that any manifold with a torsion free connection and a
non-degenerate parallel 2-form is necessarily symplectic.

We now show a few recent results from
\cite{Cahen,Rawnsley3,Rawnsley2,Rawnsley1,Vais1}, for they
constitute an important part of the theory of symplectic
connections. In order to find a smaller subspace of $\cal A$ it
was introduced in \cite{Cahen} a variational principle
\[ \int_M\, R^2\,\omega^n   \]
where
$R^2=\underline{R}_{\alpha\beta\gamma\delta}\underline{R}^{\alpha\beta\gamma\delta}$,
with $\underline{R}^{\alpha\beta\gamma\delta}=
\underline{R}_{\alpha'\beta'\gamma'\delta'}\omega^{\alpha\alpha'}
\omega^{\beta\beta'}\omega^{\gamma\gamma'}\omega^{\delta\delta'}$,
and where
\[   \underline{R}(X,Y,Z,T)=\omega(R^\nabla(X,Y)Z,T)    \]
--- a tensor in $\wedge^2T^*M\otimes S^2T^*M$.
This verifies the first Bianchi identity, because $T^\nabla=0$,
and a second Bianchi identity
\[  \cyclic_{X,Y,Z}(\nab{X}{\underline{R}})(Y,Z,T,U)=0 .    \]
The representation theory on the space of tensors like
$\underline{R}$, under the action of $Sp(2n,\R)$, has been
determined by I. Vaisman in \cite{Vais1}. It is known that the
curvature of $\nabla$ has two irreducible components. So we write
$\underline{R}=E+W$ where
\begin{eqnarray*}
\lefteqn{ E(X,Y,Z,T)=-\frac{1}{2(n+1)}\Bigl\{2\omega(X,Y)r(Z,T)+\omega(X,Z)r(Y,T)  \Bigr. }  \\
    & & \hspace{3.3cm}  \Bigl. +\omega(X,T)r(Y,Z)-\omega(Y,Z)r(X,T)-\omega(Y,T)r(X,Z)\Bigr\}.
\end{eqnarray*}
$r(X,Y)=\Tr{\left\{Z\mapsto R^\nabla(X,Z)Y\right\}}$ is the Ricci
tensor. $W$ is called the Weyl tensor and the connection is said
to be of Ricci type if $W=0$. This Weyl part of $\underline{R}$
plays a role parallel to that of the Weyl curvature tensor in
Riemannian geometry. In our case too, it is 0 in dimension 2. The
variational principle yields the field equations
\begin{eqnarray}
  \cyclic_{X,Y,Z}(\nab{X}{r})(Y,Z)=0 , \label{fieldequations}
\end{eqnarray}
having as particular solutions the Ricci type connections.

In \cite{Rawnsley1} we meet a further characterization of $\nabla
r$ and find the interesting result that, if
$(M_i,\omega_i,\nabla_i)$ are two symplectic manifolds together
with corresponding symplectic connections and such that the
symplectic $\nabla=\nabla_1+\nabla_2$ over the cartesian product
$(M_1\times M_2, \omega_1+\omega_2)$ is of Ricci type, then all
three connections must be flat.

A result proved in \cite{Rawnsley3} shows that with its standard
2-form $\C\Pro^n$ is of Ricci type --- of course the Levi-Civita
connection becomes a symplectic connection in the K\"ahlerian
framework.

Still for later purposes, we show the following example. Consider
$(\R^2,\omega)$ with the usual coordinates $z=x+iy$ and symplectic
structure $\omega=\frac{i}{2}\dx z\wedge\dx\overline{z}=\dx
x\wedge\dx y$, and let
\[ \partial_z=\frac{\partial}{\partial z}=
\frac{1}{2}\left(\frac{\partial}{\partial
x}-i\frac{\partial}{\partial y}\right),\hspace{1.1cm}
\partial_{\overline{z}}=\frac{\partial}{\partial
\overline{z}}=\frac{1}{2}\left(\frac{\partial}{\partial
x}+i\frac{\partial}{\partial y}\right). \]
\begin{prop} \label{prop2.1}
Every symplectic connection  on $(\R^2,\omega)$ is uniquely
determined by two functions $\alpha,\beta$$\in\cinf{M}(\C)$
satisfying
\[  \nab{\partial_z}{\partial_z}=\alpha\,\partial_z+\beta\,\partial_{\overline{z}}=
   \overline{\nab{\partial_{\overline{z}}}{\partial_{\overline{z}}}}  \]
and
\[ \nab{\partial_z}{\partial_{\overline{z}}}=
-\overline{\alpha}\,\partial_z-\alpha\,\partial_{\overline{z}}=
\nab{\partial_{\overline{z}}}{\partial_z}. \]
\end{prop}
The proof is elementary. Indeed, the {\it real} and torsion free
assumptions, together with
\[ \frac{i}{2}\alpha= \omega(\nab{\partial_z}{\partial_z},\partial_{\overline{z}})
  =-\omega(\partial_{z},\nab{\partial_z}{\partial_{\overline{z}}}),   \]
give us the formula.
Because sometimes is impossible to {\it keep} complex, we show the
real correspondent of the last proposition. If
\begin{eqnarray}
   \nab{\partial_x}{\partial_x}\,=\,b\,\partial_x-a\,\partial_y, \hspace{1.73cm}& &
    \ \ \  \nab{\partial_y}{\partial_y}\,=\,d\,\partial_x-c\,\partial_y,   \\
  \nab{\partial_x}{\partial_y}\,=\,c\,\partial_x-b\,\partial_y\,=\,\nab{\partial_y}{\partial_x},
  & &                                \label{eq8}
\end{eqnarray}
with $a,b,c,d:\R^2\rightarrow\R$, then
\[ \alpha=-\frac{b+d}{4}-i\frac{a+c}{4}\ \ \ \ \ \ \mbox{and}\ \ \
\ \ \ \beta=\frac{3b-d}{4}-i\frac{3c-a}{4}     . \]

\vspace{.5cm}

\subsection{Translation invariant symplectic connections}
\newcommand{\Jacob}{{\mathrm{Jac}}\,}

Here we study connections in $M=\R^{m}$. Let $s$ denote the global
frame
\[ s=\left(\frac{\partial}{\partial x_1},\ldots,\frac{\partial}{\partial x_m}\right)  \]
and let $\nabla^0=\dx$ be the trivial connection: in the example
above, $\alpha=\beta=0$.
\begin{prop} \label{prop2.2}
(i) A \co\ $\nabla$ is flat iff there exists an open cover
$\{U_i\}$ of $M$ and a collection of maps $g_i:U_i\rightarrow
GL(m,\R)$ such that
\begin{equation}                                           \label{eqn211}
 \nabla s=sg_i\dx g_i^{-1}.
\end{equation}
(ii) Let  $\sigma\in Dif\!f(\R^m)$ and
$\nabla=\sigma\cdot\nabla^0$. Then the map $g$ given by
\begin{equation}                                               \label{eqn212}
g\circ\sigma=\Jacob \sigma
\end{equation}
satisfies equation \eqref{eqn211} globally.\\
(iii) Given any  map $g:M\rightarrow GL(m,\R)$, a necessary
condition for equation \eqref{eqn212} to have solution in variable
$\sigma$, corresponding to Schwarz theorem of
mixed derivatives, is that the flat \co\ defined by formula \eqref{eqn211} is torsion free.\\
(iv) The isotropy subgroup of $\nabla^0$ is
$Dif\!f(M)_{\nabla^0}=GL(m,\R)\rtimes\R^m$.
\end{prop}
\begin{proof}
{\it (i)} The condition of $\nabla$ being flat is equivalent to
the local existence of parallel frames (solution to a system of
quasi-linear differential equations of the first
order). The result follows by straightforward computations.\\
{\it (ii)} $\nabla$  is flat since $R^\nabla=\sigma\cdot
R^{\nabla^0}=0$. Now, we have that $\sigma\cdot s=sg$ with some
$g:\R^m\rightarrow GL(m,\R)$. Then
\[ \nab{X}{sg}=\sigma\cdot\left(\nabla^0_{_{\sigma^{-1}\cdot X}}\,s\right)=0\]
for any vector field $X$. On the other hand,
\[  \nabla\,sg =(\nab{}{s})g+s\dx g=s(Ag+\dx g)  \]
if $\nabla=\nabla^0+A$. Hence $Ag=-\dx g$, which is equivalent to
$A=g\dx g^{-1}$. Now, for $\sigma=(\sigma_1,\ldots,\sigma_m)$,
since
\begin{eqnarray*}
\sigma\cdot s & = & \biggl(\dx\sigma\Bigl(\frac{\partial}{\partial
x_1} \,\mbox{\tiny{$(\sigma^{-1})$}}\Bigr),\ldots,\dx\sigma\Bigl(
\frac{\partial}{\partial x_m}\,\mbox{\tiny{$(\sigma^{-1})$}}\Bigr)\biggr)\\
 & = & \left(\frac{\partial\sigma_1}{\partial x_1}\frac{\partial}{\partial x_1}
 +\ldots+\frac{\partial\sigma_m}{\partial x_1}\frac{\partial}{\partial x_m}\,,\ldots,
 \,\frac{\partial\sigma_1}{\partial x_m}\frac{\partial}{\partial x_1}+\ldots+
 \frac{\partial\sigma_m}{\partial x_m}\frac{\partial}{\partial x_m}\right)
 \mbox{\tiny{$(\sigma^{-1})$}}\\
 & = & s\Jacob\sigma_{\mid \sigma^{-1}}\:,
\end{eqnarray*}
we may deduce formula \eqref{eqn212}.\\
{\it (iii)} Again, notice that $T^\nabla=\sigma\cdot
T^{\nabla^0}=0$ is a necessary condition for the existence of a
map $\sigma$. We just have to check this agrees with Schwarz
theorem. On one hand, if $y=\sigma(x)$,
\[\frac{\partial^2\sigma_j}{\partial x_i\partial
x_k}=\frac{\partial}{\partial
x_i}(g_{jk}\circ\sigma)=\frac{\partial g_{jk}}{\partial
y_l}\frac{\partial y_l}{\partial x_i}=\frac{\partial
g_{jk}}{\partial y_l}g_{li}.\] On the other hand, if $(e_k)$ is
the canonical basis,
\[\nab{\partial_i}{\partial_k}=-s\dx g\left(\frac{\partial}{\partial
x_i}\right)g^{-1}e_k^T=-s\frac{\partial g_{jt}}{\partial
x_i}g^{tk}.\] This is equal to $\nab{\partial_k}{\partial_i}$ iff
\[\left(\frac{\partial g_{jt}}{\partial
x_i}g^{tk}\right)g_{k\alpha}g_{i\beta}=\left(\frac{\partial
g_{jt}}{\partial x_k}g^{ti}\right)g_{k\alpha}g_{i\beta}\] or
\[ \frac{\partial g_{j\alpha}}{\partial x_i}g_{i\beta}=\frac{\partial
g_{j\beta}}{\partial x_k}g_{k\alpha} \] which is the equation we
were looking for.
\\
{\it (iv)} If $\sigma\cdot\nabla^0=\nabla^0$, then $\Jacob\sigma$
is constant by equations \eqref{eqn211} and \eqref{eqn212}.
Integrating, we find the group of affine transformations.
\end{proof}
Given a \co\ $\nabla$, notice the system of partial differential
equations $\nabla=\sigma\cdot\nabla^0$ in variable $\sigma$ is
$2^{\mathrm{nd}}$-order nonlinear. However, we have checked that
it is a $1^{\mathrm{st}}$-order linear in the entries of
$\Jacob\sigma_{\mid \sigma^{-1}}$ and easily integrated as such.
Indeed, we have just proved that solving
$\nabla=\sigma\cdot\nabla^0$ is \emph{equivalent} to solving
equation \eqref{eqn212}.

Supposing solutions $\sigma$ exist, composing them with any
translation $x\mapsto x+v,\ \,v\in\R^m$, will also give a
solution. So we may look for $\sigma$ such that $\sigma(0)=0$.
Also, assuming $g(0)=1$ is not a problem either, as one deduces
from the formula in the proposition --- it corresponds to a gauge
transformation.

We give a simple example just to illustrate the proposition:
consider the open set $\R^+\times\R$ and, in real coordinate
functions, {\it cf.} \eqref{eq8}, take the \co\ $a=c=0,\ \,d=x$
and $b=-\frac{1}{2x}$. An easy computation shows $\nabla$ is flat.
A little extra work to find the group-valued map $g$, leads then
to the problem of finding $(\sigma_1,\sigma_2)$ such that
\begin{eqnarray*}
\left[ \begin{array}{cc}
\frac{\sqrt{2\sigma_1}}{2}e^{-\frac{\sigma_2}{\sqrt{2}}} &
 -\sqrt{\sigma_1}e^{\frac{\sigma_2}{\sqrt{2}}} \\
\frac{1}{2\sqrt{\sigma_1}}e^{-\frac{\sigma_2}{\sqrt{2}}} &
\frac{\sqrt{2}}{2\sqrt{\sigma_1}}e^{\frac{\sigma_2}{\sqrt{2}}}
\end{array}\right] =
\left[ \begin{array}{cc}
\frac{\partial\sigma_1}{\partial x} & \frac{\partial\sigma_1}{\partial y} \\
\frac{\partial\sigma_2}{\partial x} &
\frac{\partial\sigma_2}{\partial y}
\end{array}\right] .
\end{eqnarray*}
Notice $\bigl\{\sigma_1,\sigma_2\bigr\}=\det\,\Jacob\sigma=1$.
This is the case where the map $g$ takes values in
$SL(2)=Sp(2,\R)$. In general, if the map is $G$-valued, then
$\nabla$ is a $G$-connection.

\vspace*{.5cm}

There is a type of connections for which we have found a solution
to the problem raised before. Consider a symplectic connection in
$\R^{2n}$ which is translation invariant, that is
$T_v\cdot\nabla=\nabla$ for all maps $T_v(x)=x+v$,
\,$v\in\R^{2n}$. Letting $\nabla=\nabla^0+A$ where $A$ is a
$\symp(2n,\R)$-valued 1-form, then we must have
\[ T_v\cdot(\nabla^0+A)=\nabla^0+T_v\cdot A=\nabla^0+A.  \]
Since $\dx T_v=\Id$, one does not take long to conclude that
$A_{x+v}=A_x$, {\it ie.} $A$ is a constant 1-form. The following
theorem has been proved with entirely different methods.
\begin{teo}[\cite{Rawnsley2}] \label{teo2.3}
Let $\nabla$ be a flat, translation invariant and symplectic \co\
on the manifold $\R^{2n}$. Suppose $\nabla=\nabla^0+A$. Then
$A(X)A(Y)=0$ for all vectors $X,Y$, and with the map
\[ \sigma(x)=x-\frac{1}{2}A(x)x   \]
we have $\nabla=\sigma\cdot\nabla^0$.
\end{teo}
\begin{proof}
First we have
\begin{eqnarray*}
 0= R^\nabla=\dx^{\nabla^0} A+A\wedge A = A\wedge A
\end{eqnarray*}
so that $[A(X),A(Y)]=0$ for any pair of vector fields. Hence, to
see $A(X)A(Y)=0$, we just have to show $A(X)A(X)=0$. Let
$X\in\R^{2n}$ be fixed and consider the 2-form
\[ \alpha(Y,Z)\,=\,\omega(A(X)Y,A(X)Z).  \]
By the torsion free assumption, $A(X)Y=A(Y)X$, so $\alpha$ also
satisfies
\begin{eqnarray*}
\alpha(Y,Z) &=&\omega(A(Y)X,A(Z)X)\\
  &=& -\omega(A(Y)A(Z)X,X)
\end{eqnarray*}
and hence, being symmetric, it must vanish ---  which implies
$A(X)A(X)=0$. This proves the first part of the theorem.

From proposition \ref{prop2.2} we have that $A=g\dx g^{-1}$ for
some global $g\in{\mathrm A}^0(Sp(2n,\R))$. Certainly, in
canonical coordinates $(x_1,\ldots,x_{2n})$
\[ A=\sum A_i\dx x_i=\dx\left(\sum x_iA_i\right)  \]
with constant $A_i$. Now let $B=\sum x_iA_i=A(x)$. Again, $\dx
B\,B=B\dx B$ so if we put
\[ g=e^{-B}=\sum_{m\geq0} \frac{(-B)^m}{m!}  \]
then $g\dx g^{-1}=\dx B=A$. According to the same proposition
\ref{prop2.2} we are left to solve the equation
\[ e^{-A(\sigma)}=\Jacob\sigma  \]
or equivalently
\[  1-A(\sigma)=\Jacob\sigma .  \]
In the canonical basis $(e_i)$ of $\R^{2n}$, this is the same as
\[ e_i-A(\sigma)e_i=\frac{\partial\sigma}{\partial x_i} .  \]
Letting $\sigma(x)=x-\frac{1}{2}A(x)x $ \,then on one hand we have
\[ A(\sigma(x))\,=\,A(x)-\frac{1}{2}A(A(x)x)\,=\,A(x)  \]
and on the other
\begin{eqnarray*}
\frac{\partial\sigma}{\partial x_i}\,=\,
e_i-\frac{1}{2}A(e_i)x-\frac{1}{2}A(x)e_i\, =\,e_i-A(x)e_i
\end{eqnarray*}
so the given map satisfies the differential equation, as we
wished.
\end{proof}

We acknowledge the help of \cite{Rawnsley2} in seeing the
$A(X)A(Y)=0$ part, in dimensions $n\geq2$. Finally, one may easily
find the set of non-zero 1-forms $A$ representing flat,
translation invariant symplectic \co s in $\R^2$, up to a scalar
factor. It is in 1-1 correspondence with the {\it non-empty} curve
\[ \Bigl\{[a:b:c:d]\in P^3(\R):\ bc-ad=0,\ b^2-ac=0\Bigr\}\backslash
\{pt\} , \] where $pt=[0:0:1:0]$, which we may compactify by
adding the trivial connection.

\vspace{1cm}

\section{Examples}

According to \cite{Obri}, a twistor space over a base space $M$ is
an almost complex manifold $Z$ together with a submersion
$f:Z\rightarrow M$ with fibres almost complex submanifolds. For
each $z$ in the fibre $Z_x=f^{-1}(x)$ we have an isomorphism
\[ \frac{T_zZ}{{\cal V}_z}\longrightarrow T_xM     \]
where ${\cal V}_z=\ker\,\dx f_z=T_zZ_x$. Then, since the vector
space $T_zZ/{\cal V}_z$ is complex, we can take this complex
structure to $T_xM$ in order to construct a map
\[ j:Z\longrightarrow {\cal J}(M). \]
Of course $f$ is a pseudo-holomorphic map with respect to some
structure on $M$ if, and only if, $j$ is constant along the
fibres.

If $(M,\omega)$ is a symplectic manifold, we shall call $Z$ an
``$\omega$-twistor space'' if the image of $j$ is in some $\zl$.
For example, given a symplectic connection $\nabla$ on $M$, the
tautology of the definition of $\jnab$ proves $\zl$ to be a true
$\omega$-twistor space.

Recall that the Siegel domain and all $G/U^l$ sit holomorphically
and separately in a Grassmannian. So we ask for an extension of
$\jnab$ to the compact $Sp(n)/U(n)$-bundle of complex, lagrangian
$n$-planes over the real symplectic $2n$-manifold $M$. Such
extension does not exist (unfortunately), although the standard
fibre is a complex symmetric space.
\begin{prop}  \label{prop3.1}
It is not possible to extend $(\zl,\jnab)$ to a bigger almost
complex manifold, of the same dimension, which is also a fibre
bundle over $M$.
\end{prop}
\begin{proof}
Assuming the extension to an almost complex space $Z$ exists, the
theory above yields a continuous map
\[ j:\overline{\zl} \longrightarrow\jm  \]
on the closure of $\zl$ in $Z$, because, if $z$ is any point on
the boundary of the twistor space, projecting to a point $x\in M$,
then $T_zZ_x$ is still a complex vector space --- it is the limit
of complex vector subspaces inside a complex vector space.

Also by continuity, we have that $\omega=\omega^{1,1}$ for $j(z)$
and the induced inner product $\omega(\ ,j(z)\ )$ has the same
signature. However, $j$ is the identity in $\zl$ so we arrive to a
contradiction.
\end{proof}

Regarding a matter of different nature, it seems that the
`non-constant' compact $\omega$-twistor spaces are not easy to
construct or describe.
\begin{prop}  \label{prop3.2}
There are no $\omega$-twistor spaces with compact fibres of
$\dim>0$ satisfying the hypothesis of theorem \ref{teo1.1} and
with the map $j$ an embbeding.
\end{prop}
\begin{proof}
Assuming $Z$ was such a space, then
\[ j:Z\longrightarrow \zl \]
would be holomorphic when restricted to each fibre. However, any
$G/U^l$ is a Stein manifold so its compact analytic submanifolds
are points.
\end{proof}

Clearly the proposition avoids the case of any holomorphic
submersion $f:Z\rightarrow M$, which induces a map $j$ constant
along the fibres.
\vspace{4mm}\\
Here we have the promised examples of twistor spaces.\\
{\bf Example 1.} Let $M=\R^2,\ \omega$ the canonical symplectic
form, $\nabla$ any symplectic connection on $M$ --- see
proposition \ref{prop2.1}, from which we use the descriptions and
notations in what follows. We want to describe $\zo_M$ in terms of
its $\db$ operator, since $\jnab$ is always integrable. There is a
simple way to see this: $R^\nabla$ is a 2-form, so it is
proportional to $\omega$. Then, since $\omega=\omega^{1,1}$ for
$j\in \zo_M$, we have $R^\nabla(j^-\ ,j^-\ )=0$, and so we apply
theorem \ref{teo1.1} to prove the claim. Otherwise, one may recall
that the Weyl part of the curvature is always zero in the two
dimensional case.

Now suppose $v\in T^{0,1}M=T^-M$ for $j$. If
$v=\frac{\partial}{\partial z}$ then $j\in -\zo={\cal Z}^1$, so we
may already assume, up to a scalar,
\[   v=\frac{\partial}{\partial\overline{z}}+w\frac{\partial}{\partial z} \]
for some $w\in\C$. The ``positive'' condition reads $
-i\omega(v,\overline{v})<0$. Since
\begin{eqnarray}
 -i\omega(v,\overline{v}) &= &\frac{1}{2}\dx z\wedge\dx\overline{z}
 \left(\frac{\partial}{\partial\overline{z}}+w\frac{\partial}{\partial z},
 \frac{\partial}{\partial z}+\overline{w}\frac{\partial}{\partial\overline{z}}\right)
 \nonumber\\
& = & \frac{1}{2}(w\overline{w}-1),              \label{disco}
\end{eqnarray}
we recover\footnote{In \cite{Tese} it is proved that the map
$J\mapsto -i$-eigenspace is holomorphic, from the Siegel domain,
with `left multiplication by $J$ on $T_JG/U^l$', to the
Grassmannian of complex $n$-planes in $\C^{2n}$} the Siegel disk
${\cal D}=\{w:\ |w|<1\}$. Because $TM$ is $\cinf{}$-trivial we
have
\[ \zo=M\times{\cal D}\stackrel{\pi}\longrightarrow M.   \]

Now working together with $T\zo\otimes\C$  let
\[ u=\frac{\partial}{\partial\overline{z}}+w\frac{\partial}{\partial z}+
{\cal P}\frac{\partial}{\partial w}+{\cal
Q}\frac{\partial}{\partial\overline{w}}   \] be a
$\jnab$-(0,1)-horizontal vector field, thus projecting to
$v=\dx\pi(u)$ and where $w$ is the fibre variable. Recall the
canonical section $\Phi\in\Gamma(\End{\pi^*TM})$ defined by
$\Phi_{j}=j$. Then
\[ \Phi v=-iv \]
where we see $v$ as a (0,1)-section of $(\pi^*TM)^c$. We can
compute the function $\mathcal{P}$ solving
\begin{eqnarray}
 \left(\pi^*\nabla_u\Phi\right)v=0.  \label{equacaohorizontal}
\end{eqnarray}
On the left hand side we have --- recall proposition \ref{prop2.1}
---
\begin{eqnarray*}
\lefteqn{\left(\pi^*\nabla_u\Phi\right)v\ =\ \pi^*\nabla_u\Phi v-\Phi\,\pi^*\nabla_u v}    \\
& = & -(i+\Phi)\pi^*\nabla_u v \\
& = &
-(i+\Phi)\left(\nabla_{\dx\pi(u)}\frac{\partial}{\partial\overline{z}}+
    u(w)\frac{\partial}{\partial z}+w\nabla_{\dx\pi(u)}\frac{\partial}{\partial z}\right)\\
& = & -(i+\Phi)\left(
\nab{\partial_{\overline{z}}}{\partial_{\overline{z}}}+
   w\nab{\partial_z}{\partial_{\overline{z}}}+\mathcal{P}\frac{\partial}{\partial z}+
   w\nab{\partial_{\overline{z}}}{\partial_z}+w^2\nab{\partial_z}{\partial_z}\right)\\
& = &
-(i+\Phi)\left(\overline{\alpha}\frac{\partial}{\partial\overline{z}}+
    \overline{\beta}\frac{\partial}{\partial z}-\overline{\alpha}w\frac{\partial}{\partial z}-
    \alpha w\frac{\partial}{\partial \overline{z}}\right.\\
& & \hspace{3cm}\left.+\mathcal{P}\frac{\partial}{\partial z}-
    \overline{\alpha}w\frac{\partial}{\partial z}-
    \alpha w\frac{\partial}{\partial \overline{z}}+w^2\alpha\frac{\partial}{\partial z}+
    w^2\beta \frac{\partial}{\partial\overline{z}}\right)\\
& = &
-(i+\Phi)\left((\overline{\beta}-2\overline{\alpha}w+\mathcal{P}+w^2\alpha)
   \frac{\partial}{\partial z}+(\overline{\alpha}-2\alpha
   w+w^2\beta)\frac{\partial}{\partial\overline{z}}\right).
\end{eqnarray*}
Equation \eqref{equacaohorizontal} says we are in the presence of
a (0,1)-vector for $j$, therefore by colinearity there exists
$\lambda\in\C$ such that
\[ (\overline{\beta}-2\overline{\alpha}w+\mathcal{P}+w^2\alpha)
   \frac{\partial}{\partial z}+(\overline{\alpha}-2\alpha w+w^2\beta)
    \frac{\partial}{\partial\overline{z}}=\lambda\left(\frac{\partial}{\partial\overline{z}}+
    w\frac{\partial}{\partial z}\right). \]
Henceforth
\[ \overline{\beta}-2\overline{\alpha}w+\mathcal{P}+w^2\alpha=\overline{\alpha}w-2\alpha w^2+
   w^3\beta  \]
and thus we get a cubic in $w$ with coefficients in
$\cinf{\Rpequeno^2}(\C)$:
\[ \mathcal{P}=-\overline{\beta}+3\overline{\alpha}w-3\alpha w^2+\beta w^3.       \]
To find the function $\mathcal{Q}$ one would have to proceed as
above but with (1,0)-vector fields.
\begin{prop}  \label{prop3.3}
(i) $f\in{\cal O}_{\zo}$ if and only if
\begin{eqnarray*}
\begin{cases}
\frac{\partial f}{\partial\overline{w}}=0\\
\frac{\partial f}{\partial\overline{z}}+w\frac{\partial
f}{\partial z}+{\cal P}(z,w)\frac{\partial f}{\partial w}=0.
\end{cases}
\end{eqnarray*}
(ii) Let $j:\R^2\rightarrow\zo$ be a section, represented in
coordinates by the map $z\mapsto(z,w(z))$. Then $j$ is
$(j,\jnab)$-holomorphic iff $w$ satisfies
\[ \frac{\partial w}{\partial\overline{z}}+w\frac{\partial w}{\partial z}-{\cal P}(z,w(z))=0. \]
\end{prop}
\begin{proof} {\it (i)}
According to the footnote, $\partial/\partial\overline{w}$ is a
(0,1)-vector field tangent
to the fibres of $\zo$, hence the first equation. The second is $u(f)=0$.\\
{\it (ii)} We consider holomorphic functions $f$ on the twistor
space, thus satisfying the system in {\it (i)}, and then claim
that $j$ is $(j,\jnab)$-holomorphic iff $f\circ j$ is holomorphic,
$\forall f$. This corresponds to
\[ \dx(f\circ j)\left(\frac{\partial}{\partial\overline{z}}+w(z)
\frac{\partial}{\partial z}\right)=0. \] Equivalently,
\begin{eqnarray*}
\frac{\partial f}{\partial\overline{z}}+\frac{\partial f}{\partial
w} \frac{\partial w}{\partial\overline{z}}+\frac{\partial
f}{\partial\overline{w}}
\frac{\partial\overline{w}}{\partial\overline{z}} +
w\frac{\partial f}{\partial z} +w\frac{\partial f}{\partial
w}\frac{\partial w}{\partial z} + w\frac{\partial
f}{\partial\overline{w}}\frac{\partial\overline{w}}{\partial z}\
=\ 0
\end{eqnarray*}
or
\[ \left(-{\cal P}(z,w)+\frac{\partial w}{\partial\overline{z}}+w
\frac{\partial w}{\partial z} \right)\frac{\partial f}{\partial
w}\ =\ 0. \] Since there exist sufficient holomorphic functions,
we are finished.
\end{proof}
\vspace{4mm} \hspace{-6mm}{\bf Remarks:} 1. By Darboux's theorem,
the proposition describes locally the twistor space of
any Riemann surface.\\
2. We give an independent proof of integrability: for the given
basis of (0,1)-vector fields, we have that $ \left[\frac{\partial
}{\partial\overline{w}},u\right]=\frac{\partial
\mathcal{P}}{\partial\overline{w}}\frac{\partial }{\partial
w}+\frac{\partial c}{\partial\overline{w}}\frac{\partial
}{\partial\overline{w}}=\frac{\partial
c}{\partial\overline{w}}\frac{\partial }{\partial\overline{w}} $
is again a (0,1)-tangent vector. (The almost complex structure
$\jnab_2=(J^h,-J^v)$ is never integrable, {\it cf.}
\cite{Tese,Sal}, so this computation confirms the correct choices
in our example.)\\
3. In the general theory of twistor spaces, a section $j$ is
$(j,\jnab)$ holomorphic if and only if it satisfies a well known
condition ({\it cf.} \cite{Rawnsley,Sal}): $ \nabla_uv\in\Gamma
T^+M,\ \,\forall u,v\in\Gamma T^+M$.
\vspace{4mm}\\
{\bf Example 1.1.} This is the trivial case; recall the \co\
$\nabla^0$ is symplectic because $M=\R^2$ is K\"ahler, so assume
$\alpha=\beta=0$. We have the following global chart:
\begin{eqnarray*}
  \zo= M\times{\cal D}\ \longrightarrow\, \ \C\times{\cal D}\hspace{7mm}\\
      \ \ (z,w) \longmapsto \ (w\overline{z}-z,w)
\end{eqnarray*}
(this map is injective if and only if $|w|\neq1$). Adding a point
at infinity on the right hand side and recalling the grassmannian
model of the twistor space, the same map composed with $1/w$ gives
a chart of ${\cal Z}^1$. Curiously, this example is the only one
for which the natural fibre chart $w$ is a {\it globally}
holomorphic function. Also, $\C\times{\cal D}$ is convex, so we
conclude $\zo_M,{\cal Z}^1_M$ with complex structure ${\cal
J}^{\nabla^0}$ are both Stein 2-manifolds.

One could also try to find the global charts for the flat torus or
cylinder.
\vspace{4mm}\\
{\bf Example 2.} This is the generalisation of example 1.1. Let
$\omega=\frac{i}{2}\sum_{k=1}^n\dx z_k\wedge\dx\overline{z}_k$. We
give a description of ${\cal Z}^0_{\Rpequeno^{2n}}$ with complex
structure arising from $\nabla^0$.

First notice that for any element $j$ we can find a basis of
$T^-M$ with vectors of the kind
\[ v_k=\frac{\partial}{\partial\overline{z}_k}+\sum_lw_{kl}\frac{\partial}{\partial z_l} ,\]
with $k=1,\ldots,n,\ w_{kl}\in\C$. For, if a linear combination of
the $\partial/\partial z_l$, only, were in $T^-M$, then the
positive condition would not be satisfied. Now, $\omega$ being
(1,1) for $j$ implies
\[ 0=\omega(v_{k_1},v_{k_2})=\frac{i}{2}(w_{k_1k_2}-w_{k_2k_1}).  \]
The positive condition is given by
\begin{eqnarray*}
0\ >\  -i\omega(v_k,\overline{v}_k) &= &\frac{1}{2}\sum_l\dx
z_l\wedge\dx\overline{z}_l
 \left(\frac{\partial}{\partial\overline{z}_k} + w_{kp}\frac{\partial}{\partial z_p},
 \frac{\partial}{\partial z_k} + \overline{w}_{kq}\frac{\partial}{\partial\overline{z}_q}
 \right) \\
& = & \frac{1}{2}\sum_l\left(-\delta_{kl}+w_{kl}\overline{w}_{kl}\right) \\
& = & \frac{1}{2}\bigl(-1 +\sum_l|w_{kl}|^2\bigr)
\end{eqnarray*}
where repeated indices in $p,q$ have denoted a sum. With respect
to the symmetric matrix $W=[w_{kl}]$ this is equivalent to
$1-WW^*>0$ and so we meet another well known description of the
Siegel domain.

Continuing to reason as in example 1 we find that a function $f$
on the twistor space is holomorphic if, and only if, \,$v_k(f)=0,\
\partial f/\partial\overline{w}_{pq}=0$. So a global chart for
$\zo_{\Rpequeno^{2n}}$ is given by the functions
\begin{eqnarray*}
 f_{pq}(z_1,\dots,z_n,w_{11},\ldots,w_{n-1,n})=w_{pq}\ ,\\
f_k(z_1,\dots,z_n,w_{11},\ldots,w_{n-1,n})=\overline{z}_kw_{kk}-z_k
\end{eqnarray*}
where $p\leq q$ and $1\leq k\leq n$.
\vspace{4mm}\\
{\bf Example 3.} Consider $M=S^2=\R^2\cup\{\infty\}$ with its
K\"ahler metric and corresponding Levi-Civita connection, which is
thus symplectic. The 2-form is $\omega=\frac{i}{2}\frac{\dx
z\wedge \dx\overline{z}}{(1+|z|^2)^2}$ so, \linebreak[4]proceeding
as in \eqref{disco}, we describe the fibres over the open set
$\R^2$ with the disk $\cal D$ again. Following the theory of
hermitian manifolds, the connection is type (1,0), {\it ie.}
transforms holomorphic sections in (1,0)-forms. Thus $\nabla$ on
$T^*M$ is determined by
\[ \nabla\dx z=\alpha\,\dx z\otimes\dx z     \]
and a conjugate version of this equation, bearing in mind $\nabla$
is real. Solving $\nabla\omega=0$ leads to
\[   \alpha=\frac{2\overline{z}}{1+|z|^2}   .    \]
Proceeding then exactly as in example 1 we find: $f\in{\cal
O}_{{\cal Z}^0_{M-\{\infty\}}}$ if and only if
\begin{eqnarray}           \label{system}
\begin{cases}
\frac{\partial f}{\partial\overline{w}}=0\\
\frac{\partial f}{\partial\overline{z}}+w\frac{\partial
f}{\partial z}+ \frac{2w(w\overline{z}-z)}{1+|z|^2}\frac{\partial
f}{\partial w}=0.
\end{cases}
\end{eqnarray}

Now let $(z_1,w_1)$ denote coordinates for the twistor space of
$M$ minus the other pole. The affine transformation on the base
$z_1=\sigma(z)=1/z$ is raised to a $\jnab$-holomorphic
transformation of the twistor space. $w_1$ is defined by requiring
that
\[ (\dx\sigma)^c\left(\frac{\partial}{\partial\overline{z}}+w
\frac{\partial}{\partial
z}\right)=\lambda\left(\frac{\partial}{\partial{\overline{z}}_1}+
w_1\frac{\partial}{\partial z_1}\right) \] for some
$\lambda\in\C$. That is, the real map $\dx\sigma$ applies a
(0,1)-$w$-vector into a (0,1)-$w_1$-vector. We find
\[ w_1=\frac{{\overline{z}}^2}{z^2}w  \]
and $(z,w)\mapsto(z_1,w_1)$ is holomorphic because one verifies by
straightforward computations that if a function $f(z_1,w_1)$
satisfies the system \eqref{system} in variables $(z_1,w_1)$ then
\[ f\left(\frac{1}{z},\frac{{\overline{z}}^2}{z^2}w \right)     \]
also satisfies the linear system in variables $(z,w)$.

We shall see in the next section that this last result is a
manifestation of $\nabla$ being $\sigma$-invariant. The latter can
be deduced by uniqueness of the Levi-Civita connection after
verifying $\sigma$ is an isometry --- which is immediate, since
$\dx z_1=-\frac{1}{z^2}\dx z$ and thus $\sigma\cdot\omega=\omega$.
\vspace{4mm}\\
{\bf Example 4.} There exist {\it compact} $\omega$-twistor
spaces: Let $\T^2$ be the torus and consider the trivial bundle
\[ Z=\T^2\times S^2\stackrel{\mathrm{pr}_1}\longrightarrow\T^2 \]
with almost complex structure $J^Z$ given by the following basis
of (0,1)-tangents: the vectors
\[ \frac{\partial}{\partial\overline{z}}+\frac{|t|}{1+|t|^2}
\frac{\partial}{\partial z}\mbox{\ \ \ and\ \ \ }
\frac{\partial}{\partial\overline{t}}.
\] $z$ is the usual chart of $\R^2$ and $t$ is a fixed affine coordinate of
$S^2=\Pro^1(\C)$. Note that, for $t\neq 0$, we have
\[ \frac{|\tfrac{1}{t}|}{1+|\tfrac{1}{t}|^2}=\frac{|t|}{1+|t|^2}<1 ,\]
so $J^Z$ is well defined and preserves the natural splitting of
$TZ$. Moreover, it is compatible with the canonical symplectic
structure of $Z$. Notice $J^Z$ is not integrable, but this is not
important for our purposes.

Hence $Z$ is an $\omega$-twistor space. The map
$j:Z\rightarrow{\cal Z}^0_{\Tpequeno^2}=\T^2\times {\cal D}$
induced by $\dx\mathrm{pr}_1$ and the \C-vector bundle
$TZ/\ker\dx\mathrm{pr}_1$ identifies with
\[ j(z,t)=\left(z,\frac{|t|}{1+|t|^2}\right) \]
by construction. For the reader to compare with proposition
\ref{prop3.2}, note that $j$ is not even open along the fibers.

\vspace{1cm}

\section{A holomorphic map}

Let $(M,\omega),(M_1,\omega_1)$ be two symplectic manifolds and
$\sigma:M\rightarrow M_1$ a symplectomorphism. Then $\sigma$
induces an invertible transformation from $\zl_{_M}$ onto
$\zl_{_{M_1}}$ preserving the fibres, {\it ie.} a map $\Sigma$
such that the diagram
\begin{eqnarray*}
\begin{array}{ccc}
\zl_{_M} & \stackrel{\Sigma}\longrightarrow & \zl_{_{M_1}} \\
\pi\downarrow \hspace{3mm}  & & \hspace{3mm} \downarrow\pi_1 \\
M & \stackrel{\sigma}\longrightarrow & M_1
\end{array}
\end{eqnarray*}
commutes. Indeed, for any $y\in M_1,\
j\in\pi^{-1}(\sigma^{-1}(y))$ we define
\[ \Sigma(j)=\dx\sigma\circ j\circ\dx\sigma^{-1}  \]
an element in $\pi_1^{-1}(y)$. It is trivial to check $\Sigma$ is
well defined and invertible.

\newcommand{\jnabum}{{\cal J}^{\nabla^1}}       
\newcommand{\jnabdois}{{\cal J}^{\nabla^2}}       

Assume $\zl_{_M},\zl_{_{M_1}}$ have twistor almost complex
structures $\jnab$ and $\jnabum$, respectively, where
$\nabla^1=\sigma\cdot\nabla$ and $\nabla$ is a given symplectic
connection. We then have the following result.
\begin{teo}  \label{teo4.1}
$\Sigma$ is pseudo-holomorphic.
\end{teo}
\begin{proof}
Notice that $\Sigma$, when restricted to each fibre, extends to a
linear map between $\End{T_{\sigma^{-1}(y)}M}$ and $\End{T_yM_1}$.
Hence
\begin{eqnarray*}
 \dx\Sigma(jA) &=& \Sigma(jA) \\
  &=& \Sigma(j)\Sigma(A)\ =\  \Sigma(j)\,\dx\Sigma(A)
\end{eqnarray*}
and we may conclude the map is {\it vertically}
pseudo-holomorphic.

Now suppose $\Sigma_*{\cal H}^{\nabla}={\cal H}^{\nabla^1}$. Using
the isomorphism $\dx\pi_1:{\cal
  H}^{\nabla^1}\rightarrow\pi_1^*TM_1$, we have
\begin{eqnarray*}
\dx\pi_1\circ\dx\Sigma\;J^h_j & = & \dx\sigma\circ\dx\pi\left((\dx\pi)^{-1}\,j\,\dx\pi\right)\\
 & = & \dx\sigma\circ j\circ\left(\dx\sigma^{-1}\circ\dx\sigma\right)\circ\dx\pi  \\
 & = & \Sigma(j)\;\dx(\sigma\circ\pi) \\
 & = & \Sigma(j)\;\dx\pi_1\circ\dx\Sigma \\
 & = & \dx\pi_1\;J^h_{\Sigma(j)}\dx\Sigma.
\end{eqnarray*}
So the theorem follows after we prove $\Sigma_*{\cal
H}^{\nabla}={\cal H}^{\nabla^1}$, which is exactly the case when
we consider the particular connection $\nabla^1$.

Fix a real symplectic vector space $V$ and let $F,F_1$ be,
respectively, the frame bundles of $M$ and $M_1$. Consider the
$G$-equivariant map
\begin{eqnarray*}
 \Lambda & : & F\longrightarrow F_1   \\
   & & p\longmapsto\dx\sigma\circ p
\end{eqnarray*}
where the points $p:V\rightarrow T_xM$ are linear isomorphisms. If
$s:U\rightarrow F$ is a section on a neighborhood $U$ of $x\in M$,
then
\[ s_1=\Lambda\circ s\circ \sigma^{-1}:\sigma(U)\longrightarrow F_1 \]
is a section on a neighborhood of $\sigma(x)$. We wish to show
first that $\Lambda$ preserves the horizontal distributions
induced by the connections. Let $\alpha,\alpha_1$ denote the
connection 1-forms on $F$ and $F_1$.
\[ \nab{X_x}{s}=s(s^*\alpha)(X_x)  \]
and
\begin{eqnarray*}
 (\sigma\cdot\nabla)_{_{Y_{\sigma(x)}}}s_1 & = &   s_1(s_1^*\alpha_1)(Y_{\sigma(x)})\\
    & = & \Lambda\circ s\circ \sigma^{-1}_{\sigma(x)}\left[(\Lambda
    \circ s\circ \sigma^{-1})^*\alpha_1\right] (Y_{\sigma(x)}) \\
 & = & \dx\sigma\, s(s^*\Lambda^*\alpha_1)\,\dx\sigma^{-1}(Y_{\sigma(x)}) \\
 & = & \dx\sigma\, s(s^*\Lambda^*\alpha_1)(\sigma^{-1}\cdot Y)_x.
\end{eqnarray*}
On the other hand, since $(\sigma^{-1}\cdot
s_1)_x=\dx\sigma^{-1}({s_1}_{\sigma(x)})=s_x$, we have
\begin{eqnarray*}
 (\sigma\cdot\nabla)_{_{Y_{\sigma(x)}}}s_1 & = &
    \sigma\cdot\left(\nabla_{\sigma^{-1}\cdot Y} \,\sigma^{-1}\cdot s_1\right)_{\sigma(x)} \\
 & = &  \dx\sigma\left(\nabla_{(\sigma^{-1}\cdot Y)_x}\, s\right)  \\
 & = &  \dx\sigma\, s(s^*\alpha)(\sigma^{-1}\cdot Y)_x.
\end{eqnarray*}
Henceforth $s^*\Lambda^*\alpha_1=s^*\alpha$ and we prove the claim
that $\ker\alpha_1=\Lambda_*\ker\alpha$ by taking horizontal
frames along paths in $M$ passing through $x$. (With vertical
fundamental vector fields one can actually see further that
$\Lambda^*\alpha_1=\alpha$.)

Finally let $\zeta:F\rightarrow Z$ be the once introduced fibre
bundle ({\it cf.} first section, formula (\ref{pium})) with bundle
map
\[ \zeta(p)=pJ_0p^{-1},  \]
where $J_0\in G/U^l$ is some compatible complex structure of $V$.
Clearly
\[ \Sigma\circ\zeta(p)=\dx\sigma\,pJ_0p^{-1}\dx\sigma^{-1}=\zeta_1\circ\Lambda(p) \]
and we know the $\zeta$ preserve the horizontal tangent bundles:
\[ \zeta_*\ker\alpha={\cal H}^\nabla, \hspace{1cm}
   {\zeta_1}_* \ker\alpha_1={\cal H}^{\sigma\cdot\nabla}.       \]
Now it is no longer difficult to see that $\Sigma_*{\cal
H}^\nabla={\cal H}^{\sigma\cdot\nabla}$.
\end{proof}
\vspace{2mm}

We notice that the construction and results above are true for the
general twistor space $\jm$. Indeed, the proof does not mention
any particular feature of symplectic manifolds.
\vspace{.4cm}\\
\textbf{Remark.} An application of the last theorem is the result
at the end of example 3 in section 3. The theorem confirms that
the PDE system given there is preserved under the change of affine
coordinates in $S^2$. It also applys in the following strictly
real situation: since $(\R^2,\omega)$ is symplectomorphic to the
Poincar\'e disk $({\cal D},\omega_1)$, where
\[ \omega_1=\frac{i}{2}\frac{\dx z\wedge\dx\overline{z}}{(1-|z|^2)^2} , \]
we can study $\zl_{\cal D}$ using the theorem and example 1 in
section 3 (it corresponds to find the Darboux coordinates in $\cal
D$ and the respective connection's parameters). \vspace{.4cm}

There is a partial converse to the theorem, which is only valid in
the symplectic category. In the following we assume all the
previous setting.
\begin{coro} \label{coro4.1}
Let $\nabla^2$ be any symplectic connection on $M_1$ and suppose
$\Sigma:(\zl_M,\jnab)$\linebreak$\rightarrow(\zl_{M_1},\jnabdois)$
is holomorphic. Then
\[ \nabla^2=\sigma\cdot\nabla  \]
ie\, $\nabla^2$ is in the affine transformation orbit of $\nabla$.
\end{coro}
\begin{proof}
We have
\[  \jnabdois=\dx\Sigma\circ\jnab\circ\dx\Sigma^{-1}={\cal J}^{\sigma\cdot\nabla} \]
so the result follows by theorem \ref{teo1.2}.
\end{proof}
We remark that the theorem has the apparent merit of transforming
a $2^{\mathrm{nd}}$ order PDE's problem into a $1^{\mathrm{st}}$
order one.

\vspace{1cm}

\section{The metric}

In order to introduce a Riemannian structure on the twistor space
$\zo_M$ we need a further amount of theory from \cite{Obri}.
Recall the exact sequence \eqref{eq3}, where $E=\pi^*TM$ is a
vector bundle over $\zo_M$ with canonical complex structure
$\Phi$. Also important to recall here are propositions
\ref{prop1.1} and \ref{prop1.2}.

Let $\nabla$ be a symplectic linear connection on the given
$2n$-dimensional symplectic manifold $M$. Let $P\in{\mathrm
A}^1({\cal V})$ denote the projection with kernel ${\cal
H}^\nabla$, induced by the \co. Via the identity
\[ {\cal V}_j=\left\{ A\in\symp(E_j,\pi^{-1}\omega):\ \,A\Phi_j=-\Phi_jA\right\} \]
$P$ can also be seen as an endomorphism-valued 1-form on the
twistor space. We may thus define a new connection on $E$ by
\[ D=\pi^*\nabla-P ,  \]
which turns $\pi^*\nabla\Phi=[P,\Phi]$ equivalent to
\[ D\Phi=0 . \]
It follows that $D$ on $\End{E}$ preserves ${\cal V}$ and hence
$DJ^v=0$. Indeed, this connection is symplectic because its
difference to an obviously symplectic connection $\pi^*\nabla$
stays within $\symp(E,\pi^{-1}\omega)$, and hence, as a
derivation, acts trivially on the 2-form.

The isomorphism $\pi_*:{\cal H}^\nabla\rightarrow E$ allows us to
transfer $D$, to give rise to a new \co\ $D$ on ${\cal H}^\nabla$
satisfying
\begin{eqnarray*}
(DJ^h)X & = & \pi_*^{-1}\left(D(\pi_* J^hX)\right)-J^h\pi_*^{-1}\left(D\pi_*X\right) \\
 & = & \pi_*^{-1}\left(D\Phi\right)\pi_*X \ =\ 0 .
\end{eqnarray*}
Henceforth we have defined a \C-linear connection on $T\zo={\cal
V}\oplus{\cal H}^\nabla$ preserving this splitting, exactly in the
same lines of the general twistor theory (\cite{Obri}). Since
$\pi_*$ resulted in a parallel and \C-linear isomorphism, one
often identifies ${\cal H}^\nabla$ with $E$.

\newcommand{\DD}[2]{{D_{_{#1}}#2}}               

Now we need the following theorem valid in general in ${\cal
J}(M)$ and which we may improve in a little detail.
\begin{teo}[\cite{Obri}]   \label{teo5.1}
The connection $D$ on the tangent bundle of $\zo_M$ has torsion
whose vertical part is the projection of
$\pi^*R^\nabla-\frac{1}{2}\,P\wedge P$ into $\cal V$, and whose
horizontal part is $\pi^*T^\nabla-P\wedge\dx\pi$ after identifying
$\hnab$ with $E$.
\end{teo}
Since $[\m_{J},\m_{J}]\subset\gl(2n,J)$, \,{\it cf.} section 1,
one concludes that the vertical part of $T^D$ is just
$P(\pi^*R^\nabla)$. Also notice we are already assuming $\nabla$
is torsion free, so both formulas in the theorem can be
simplified.

The present section is devoted to the study of a natural
Riemannian structure on $\zo_M$, whose analogous construction in
`Riemannian twistor theory' has been already considered in
\cite{Rawnsley}. To see which twistor spaces of that kind over a
4-manifold admit a K\"ahler metric one may consult \cite{Hit}. For
the symplectic case, especially $\rdoisn$ canonical, one may also
consult \cite{Nannicini1}.

Recall that $G/U^0$ is a Hermitian symmetric space, hence
K\"ahlerian. With the help of the Killing form and a Cartan's
decomposition of $\symp(2n,\R)=\uni^{0}\oplus\m_{J}$ one defines a
symplectic form on $\zo$ by
\[ \Omega^\nabla\, =\, t\,\pi^*\omega-\tau ,\]
where $t\in]0,+\infty[$ is fixed and
\[ \tau(X,Y)\ =\ \frac{1}{2}\Tr{(PX)\Phi(PY)} . \]
The following is trivial to check.
\begin{lema}   \label{lema5.1}
$\Omega^\nabla$ is non-degenerate and $\jnab$ is compatible with
it. The induced metric is positive definite.
\end{lema}
Although the parameter $t$ will not teach us anything special
about the twistor space, besides that it could also give a
pseudo-metric, it may become important at some moment.
\begin{prop}  \label{prop5.1}
For any $X,Y,Z\in T\zo$
\[  \dx\tau(X,Y,Z)\ =\
-\frac{1}{4}\Tr{}\left(R^{\pi^*\nabla}_{_{X,Y}}\circ\pi^*\nab{Z}{\Phi}+
R^{\pi^*\nabla}_{_{Y,Z}}\circ\pi^*\nab{X}{\Phi}+R^{\pi^*\nabla}_{_{Z,X}}\circ
\pi^*\nab{Y}{\Phi}\right) .\]
\end{prop}
\begin{proof}
Let us first see $D\Omega^\nabla=0$. By previous remarks we are
left to check $D\tau=0$.
\begin{eqnarray*}
D_{_X}\tau\,(Y,Z) & = & X(\tau(Y,Z))-\tau(\DD{X}{Y},Z)-\tau(Y,\DD{X}{Z})\\
  & = & X(\tau(Y,Z))-\tfrac{1}{2}\Tr{\left(P(\DD{X}{Y})\Phi PZ+PY\Phi P(\DD{X}{Z})\right)}\\
 & = & \,X(\tau(Y,Z))-\tfrac{1}{2}\Tr{\,\DD{X}{(PY\,\Phi\,PZ)}}\\
 & = &
  X(\tau(Y,Z))-\dx\left(\tfrac{1}{2}\Tr{(PY\,\Phi\,PZ)}\right)(X)\ =\ 0.
\end{eqnarray*}
Now, it is well known that
\[  \dx\tau(X,Y,Z)  =
\tau(T^D_{_{X,Y}},Z)+\tau(T^D_{_{Y,Z}},X)+\tau(T^D_{_{Z,X}},Y) .
\] Since
\begin{eqnarray*}
 \tau(T^D_{_{X,Y}},Z) & = &
 \dfrac{1}{4}\Tr{}\left([PT^D_{_{X,Y}},\Phi]PZ\right)\\
&=& -\dfrac{1}{4}\Tr{}\left(\pi^*R^\nabla_{_{\pi_*X,\pi_*Y}}[PZ,\Phi]\right)\\
&=&
-\dfrac{1}{4}\Tr{\left(R^{\pi^*\nabla}_{_{X,Y}}\circ\pi^*\nab{Z}{\Phi}\right)}
,
\end{eqnarray*}
the result follows.
\end{proof}
\begin{teo}  \label{teo5.2}
$\Omega^\nabla$ is closed if and only if $\nabla$ is flat. In such
case, $\zo_M$ is a K\"ahler manifold.
\end{teo}
\begin{proof}
Since $\dx\pi^*\omega=0$, we only have to do an analysis of
$\dx\tau$ on four cases --- with three horizontal or vertical
tangent vectors $X,Y,Z$.

The only possible non-trivial case is say $X,Y$ horizontal and $Z$
vertical. Then, since $\tau$ on ${\cal V}$ is non-degenerate,
$\dx\tau(X,Y,Z)= \tau(T^D_{_{X,Y}},Z)=0$ for all those $X,Y,Z$ iff
$P(T^D)=0$. Equivalently, $[\pi^*R^\nabla,\Phi]=0$, or
\[ [R^\nabla_x,j]=0,\hspace{7mm} \forall j\in\pi^{-1}(x),\ \,x\in M .  \]
Now, for any $J$ compatible with $(\R^{2n},\omega)$, let
$\uni^0_{J}$ be the unitary Lie algebra
$\symp(2n,\R)\cap\gl(2n,J)$. It is then trivial to see that
\[\h=\bigcap_{J\in G/U^0}\ \uni^0_{J} \]
is a $G$-module under the adjoint action. Because $\symp(2n,\R)$
is irreducible, we have $\h=0$ and thus the `only if' part of the
theorem.

For the last part of the theorem we recall that $R^\nabla=0$
implies integrability of the almost complex structure $\jnab$ as
well.
\end{proof}

Notice $D$ is always Hermitian, $\Omega^\nabla$ may be
K\"ahlerian, but $T^D$ is never 0. Thus the (0,1) part of $D$
cannot be the $\db$ operator.

Let $\langle\ ,\ \rangle$ be the induced metric, so that
\[ \langle X ,Y \rangle\ =\ t\,\pi^*\omega(X,\jnab Y)+\frac{1}{2}\Tr{(PXPY)} \]
and thus ${\cal H}^\nabla\perp{\cal V}$. Let $\cdot^v$ denote the
vertical part of any tangent-valued tensor.
\begin{teo}  \label{teo5.3}
(i) The Levi-Civita \co\ of $\langle\ ,\ \rangle$ is given by
\[ \DCON_{_X}Y\ =\ \DD{X}{Y}-PY(\pi_*X)-\frac{1}{2}\pi^*R^{\ v}_{_{X,Y}}+S(X,Y) \]
where $S$ is symmetric and defined both by
\[  \langle S^v(X,Y),A\rangle\ =\ \langle A\pi_*X,\pi_*Y\rangle,\ \ \ \ \ \ \
      \forall A\in {\cal V} , \]
and
\[  \langle S^h(X,B),Y\rangle\ =\ \frac{1}{2}\langle\pi^*R^{\
v}_{_{X,Y}},B\rangle,\ \ \ \ \ \ \ \forall Y\in{\cal H}^\nabla .
\] Hence for $X,Y\in{\cal H}^\nabla$ and $A,B\in{\cal V}$ we have
\begin{eqnarray*}
S^v(X,A)=S^v(A,B)=0 , \\
S^h(X,Y)=S^h(A,B)=0 .
\end{eqnarray*}
(ii) The fibres $\pi^{-1}(x),\ \, x\in M$, are totally geodesic in $\zo_M$.\\
(iii) If $\nabla$ is flat, then $\DCON\jnab=0$.
\end{teo}
\begin{proof}
{\it (i)} Note that $S^h$ is symmetric by definition and that, to
confirm $S^v$ is symmetric, we just have to check every $A\in{\cal
V}$ is self-adjoint:
\begin{eqnarray*}
\langle A\pi_*X,\pi_*Y\rangle & = & t\,\omega(A\pi_*X,\Phi\pi_* Y) \\
  & = & t\,\omega(\pi_*X,\Phi A\pi_*Y)\ =\ \langle\pi_*
  X,A\pi_*Y\rangle .
\end{eqnarray*}
Now let us see the torsion condition:
\begin{eqnarray*}
T^{\DCON}(X,Y) & = & T^D(X,Y)-PY(\pi_*X)-\frac{1}{2}\pi^*R^{\
   v}_{_{X,Y}}+S(X,Y)+\\
& & \hspace{2.5cm}+PX(\pi_*Y)+\frac{1}{2}\pi^*R^{\ v}_{_{Y,X}}-S(Y,X)\\
   & = & T^D(X,Y)+P\wedge\dx\pi(X,Y)-\pi^*R^{\ v}_{_{X,Y}}\ =\ 0 .
\end{eqnarray*}
For the metric condition it is wise, from now on, to let $X,Y,Z$
denote horizontal and $A,B,C$ vertical vector fields. We already
know $D$ is Hermitian, so to simplify computations let
$\xi=\DCON-D$. Then
\begin{eqnarray*}
\xi_{_X}Y\, =\, -\frac{1}{2}\pi^*R^{\
v}_{_{X,Y}}+S^v(X,Y),\hspace{1.8cm}\xi_{_X}A\, =\, -AX+S^h(X,A) , \\
\xi_{_A}X\, =\, S^h(X,A),\hspace{1.8cm}\xi_{_A}B\, =\, 0
\hspace{2.5cm}
\end{eqnarray*}
and thus in particular, from the last formula, we deduce {\it
(ii)}. Now
\[ \DCON_{_X}\langle\ ,\ \rangle(Y,Z)\ =\
-\langle\xi_{_X}Y,Z\rangle-\langle Y,\xi_{_X}Z \rangle\ =\ 0 , \]
\begin{eqnarray*}
\DCON_{_X}\langle\ ,\ \rangle(Y,A) & = & -\langle\xi_{_X}Y,A\rangle-\langle Y,\xi_{_X}A \rangle  \\
 & = &  \frac{1}{2}\langle\pi^*R^{\ v}_{_{X,Y}},A\rangle-\langle
 S^v(X,Y),A\rangle +\langle Y,AX\rangle-\langle Y,S^h(X,A) \rangle\ =\
 0 ,
\end{eqnarray*}
\begin{eqnarray*}
-\DCON_{_X}\langle\ ,\ \rangle(A,B)\ =\   \langle\xi_{_X}A,B
 \rangle+\langle A,\xi_{_X}B\rangle\ =\ 0 ,
\end{eqnarray*}
\begin{eqnarray*}
 -\DCON_{_A}\langle\ ,\ \rangle(X,Y) & = & \langle S^h(X,A),Y \rangle+\langle X,S^h(Y,A)\rangle \\
  & = & \frac{1}{2}\langle\pi^*R^{\ v}_{_{X,Y}},A \rangle+
  \frac{1}{2}\langle\pi^*R^{\ v}_{_{Y,X}} ,A\rangle\ =\ 0 ,
\end{eqnarray*}
\begin{eqnarray*}
 -\DCON_{_A}\langle\ ,\ \rangle(X,B)\ =\ \langle\xi_{_A}X,B
  \rangle+\langle X,\xi_{_A}B\rangle\ =\ 0 ,
\end{eqnarray*}
and finally
\[ -\DCON_{_A}\langle\ ,\ \rangle(B,C)\ =\ 0 . \]
{\it (iii)} We already know this, but we are glad to confirm: if
$\nabla$ is flat then $S^h=0$. Hence for {\it all} vector fields
\[ \DCON_{_X}\jnab Y\ =\ \jnab\DD{X}{Y}-\jnab PY(\pi_*X)+S^v(X,\jnab
Y) . \] It is an easy task to show $\langle S^v(X,\jnab Y),A
\rangle=\langle\jnab S^v(X,Y),A\rangle$ (an identity also
\linebreak[3]following from the theory of the $2^{\mathrm{nd}}$
fundamental form in K\"ahler geometry), so we are finished.
\end{proof}
One may write $S^v$ explicitly and construct a
symplectic-orthonormal basis of ${\cal V}$ induced by a given such
basis on ${\cal H}^\nabla$. We show the first of these assertions.
\begin{prop}  \label{prop5.2}
For $X,Y$ horizontal
\[ S^v_j(X,Y)\,\ =\ -\frac{t}{2}\Bigl\{\omega(X,\,\
)jY+\omega(jY,\,\ )X+\omega(jX,\,\ )Y+\omega(Y,\,\ )jX\Bigr\}  .\]
\end{prop}
\begin{proof}
Since this formula is clearly symmetric we just have to verify
that $S^v(X,X)\in{\cal V}$ and $\langle
S^v(X,X),A\rangle$\linebreak[3]$=\langle AX,X\rangle$ for any
vertical vector $A\in{\cal V}$. For the first part
\begin{eqnarray*}
\omega(S^v_{X,X}Y,Z) &=&-t\omega\bigl(\omega(X,Y)jX+\omega(jX,Y)X,Z\bigr) \\
& = &
-t\bigl\{\omega(X,Y)\omega(jX,Z)+\omega(jX,Y)\omega(X,Z)\bigr\}\
=\ \omega(S^v_{X,X}Z,Y)
 \end{eqnarray*}
and
\begin{eqnarray*}
S^v_{X,X}j\, &=& -t\bigl\{\omega(X,j\ \: )jX+\omega(jX,j\ \: )X\bigr\}\\
&=& tj\bigl\{\omega(jX,\ )X+\omega(X,\ )jX\bigr\}\ =\
-jS^v_{X,X}\,.
\end{eqnarray*}
Now let $(e_1,\ldots,e_n,je_1,\ldots,je_n)$ be an orthonormal and
symplectic basis of $\hnab_j\simeq T_{\pi(j)}M$. Then
\begin{eqnarray*}
\langle S^v(X,X),A\rangle &=& \frac{1}{2}\Tr{S^v(X,X)A}\\
&=&\dfrac{1}{2}\omega\bigl(S^v(X,X)Ae_i,je_i\bigr)+
                   \dfrac{1}{2}\omega\bigl(S^v(X,X)Aje_i,j^2e_i\bigr)\\
&=&\omega\bigl(S^v(X,X)Ae_i,je_i\bigr)\\
&=&-t\omega\bigl(\omega(X,Ae_i)jX+\omega(jX,Ae_i)X,je_i\bigr)\\
&=& t\omega(AX,\omega(jX,je_i)e_i - t\omega(AX,\omega(jX,e_i)je_i)\\
&=& t\omega(AX,jX)\ =\ \langle AX,X\rangle .
\end{eqnarray*}
\end{proof}

\vspace{.5cm}

\subsection{K\"ahlerian twistor spaces}

Next we present a result about the sectional curvature of the
K\"ahlerian twistor space $\zo$. Since the result is not used
anymore we do not show its long proof. Until the end of the
subsection assume $R^\nabla=0$.
\vspace{3mm}\\
\textit{Theorem. Let $\Pi$ be a 2-plane in $T_j\zo$ spanned by the
orthonormal basis $\{X+A,Y+B\},\ X,Y\in{\cal H}^\nabla,\
A,B\in{\cal V}$. Then the sectional curvature of $\Pi$ is
\begin{eqnarray*}
k_j(\Pi) & = & -\langle R^{\DCON}(X+A,Y+B)(X+A),Y+B\rangle \\
  & = & \frac{1}{2}\Big(\Vert X\Vert^2\Vert
  Y\Vert^2+3t^2\omega(X,Y)^2-\langle X,Y\rangle^2\Big)+\\
  & & \hspace{8mm} +\Vert BX-AY\Vert^2+2\langle[A,B]X,Y\rangle-\Vert [A,B]\Vert^2
\end{eqnarray*}
where $[\,,\,]$ is the commutator bracket. Thus
\[
 k_j(\Pi)\ \biggl\{\begin{array}{ll} >0 & \mbox{ for }\Pi\subset{\cal
 H}^\nabla\\  <0 & \mbox{ for }\Pi\subset{\cal V}\ . \end{array} \biggr.
\]}

We remark that the second part of the theorem can be obtained
immediately from Gauss-Codazzi's equations. First, notice that the
horizontal distribution is integrable when $\nabla$ is flat. Then
the horizontal leaves are immediately seen to have $\pi^*\nabla$
for Levi-Civita \co\  with the induced metric, and hence they are
flat. Finally, for $X,Y$ horizontal and orthonormal, and being $S$
the $2^{\mathrm{nd}}$ fundamental form, a formula of Gauss says
\begin{eqnarray*}
 k_j\{X,Y\} & = &\Vert S(X,Y)\Vert^2-\langle S(X,X),S(Y,Y) \rangle \\
 & = & \langle S(X,Y)X,Y\rangle-\langle S(X,X)Y,Y\rangle\ =\ \mbox{etc}
\end{eqnarray*}
which is positive, as we may deduce following proposition
\ref{prop5.2}. For the totally geodesic vertical fibres of $\zo$,
we recall that $-\Vert[A,B]\Vert^2$ is the sectional curvature of
the hyperbolic space $Sp(2n,\R)/U(n)$.

One can find the Cauchy-Riemann operator on the tangent bundle of
$\zo$. We shall proceed to do this, hoping to bring further
understanding to the K\"ahlerian case.
\begin{prop}   \label{prop5.3}
(i) A tangent vector field $Y$ on $\zo_M$ is holomorphic iff
\[ \DD{X}{Y}+\jnab\DD{\jnab X}{Y}-2(PY)\pi_*X\ =\ 0,\ \ \ \ \ \
\forall X\ . \]
(ii) ${\cal H}^\nabla$ is a holomorphic subvector bundle of $T\zo$.\\
(iii) $R^D$ is a (1,1)-form.
\end{prop}
\begin{proof}
{\it (i)} It is well known that $\db=\,''\circ\DCON$ when we see
the tangent space as a $\C$-vector bundle. Hence
\begin{eqnarray*}
\db_{_{X+i\jnab X}}(Y-i\jnab Y) & = & \DCON_{_X}Y+\DCON_{_{\jnab
X}}\jnab Y+i\left(
\DCON_{_{\jnab X}}Y-\jnab \DCON_{_X}Y\right) \\
& = & \DCON_{_X}Y+\jnab \DCON_{_{\jnab
X}}Y-i\jnab\left(\DCON_{_X}Y+\jnab \DCON_{_{\jnab X}}Y\right)\ .
\end{eqnarray*}
Therefore $\db$ operates as the real part of the above, which is
equal to
\begin{eqnarray*}
 D_{_X}Y-(PY)\pi_*X+S(X,Y)+\jnab D_{_{\jnab X}}Y
-\jnab(PY)\pi_*\jnab X+\jnab S(\jnab X,Y) \\
=\   D_{_X}Y+\jnab D_{_{\jnab X}}Y-2(PY)\pi_*X\ .
\end{eqnarray*}
{\it (ii)} We have seen $D$ is a Hermitian connection on ${\cal
H}^\nabla\simeq E$. From the formula above we immediately find
that $D$ determines a $\db$-operator on $E$ coinciding with $\db$,
hence integrable. Moreover, by a famous
theorem of Koszul-Malgrange (\cite{Koba}), $R^D$ must not have (0,2)-part.\\
{\it (iii)} This follows from {\it (ii)}. However, one may argue
as in corollary \ref{coro5.1}, formula \eqref{eqpi}.
\end{proof}
Notice ${\cal V}\subset\End{E}$ also inherits an integrable
complex structure as a manifold, by part {\it (ii)}. However, this
has no longer anything to do with $\DCON$ or $\jnab$.

In conclusion, the K\"ahlerian twistor space $\zo_M$ has
holomorphic charts in $\C^n\times\C^{\frac{1}{2}n(n+1)}$ like
\[ H\times W\hspace{1cm}\mbox{or}\hspace{1cm}U\times V \]
with $H\times\{w\}$ horizontal and $\{x\}\times V$ vertical, but
never a chart of the kind $H\times V$. This is not new though; it
agrees with the fact that the bundle projection $\pi$ is never
holomorphic.

\vspace{.5cm}

\subsection{Twistor space of a Riemann surface}

Until the end of this section assume $(M,\omega,J_0)$ is a Riemann
surface. Then there are various ways to describe $\zo_M$. For
example, combining the well known isomorphism \linebreak[4]
$\,j\mapsto(j+J_0)^{-1}(j-J_0)$, valid in general, with an extra
property of real dimension 2, one may deduce easily that $\zo_M$
is diffeomorphic to the radius 1 disk bundle of $T^+M\otimes_c
T^+M$ (see \cite{Tese}). This transformation is particularly
suitable for the study of $\jnab$ with $\nabla$ reducible to
$U(1)$: we then get a biholomorphism between $\omega$-twistor
spaces.

Suppose $M$ is connected, orientable and compact. Then its Euler
characteristic is equal to $2-2g$ where $g$ is the genus of $M$.
We know a way to embed $\zo_M$ in $\Pro^1(TM\otimes\C)$. Since
this is associated to an even Euler number bundle, we may use a
result from \cite{Dusa} on the classification of sphere bundles
over Riemann surfaces to conclude that it is diffeomorphic to the
trivial bundle $M\times S^2$. Hence it yields that $M$
parameterizes a disc flowing inside $S^2$, the twistor's fibres,
with boundary {\it the} principal $U(1)$-bundle of frames.

Here is a corollary of theorem \ref{teo5.1} concerning the complex
structure of twistor space.
\begin{coro}  \label{coro5.1}
If $M$ is a Riemann surface, then ${\cal H}^\nabla$ and ${\cal V}$
are holomorphic line bundles over $\zo_{_M}$.
\end{coro}
\begin{proof}
Let $D=\pi^*\nabla-P$ be the connection defined earlier, induced
here by the Levi-Civita connection $\nabla$ of $\omega(\ ,J_0\ )$.
First we compute in any dimension
\begin{eqnarray*}
\dx^{\pi^*\nabla}P(X,Y) & = & \pi^*\nab{X}{(PY)}-\pi^*\nab{Y}{(PX)}-P[X,Y]\\
 &= & \pi^*\nab{X}{PY}-\pi^*\nab{Y}{PX}+P\left(T^D(X,Y)-\DD{X}{Y}+\DD{Y}{X}\right)\\
&= & \pi^*\nab{X}{PY}-\pi^*\nab{Y}{PX}+P(\pi^*R^\nabla_{_{X,Y}}) \\
& & \hspace{1.3cm} -\pi^*\nab{X}{PY}+[PX,PY]+\pi^*\nab{Y}{PX}-[PY,PX]\\
  &=& P(\pi^*R_{_{X,Y}})+2[PX,PY].
\end{eqnarray*}
Hence, by a well known formula on the curvature, we have
\begin{eqnarray}
R^D & =& R^{\pi^*\nabla}-\dx^{\pi^*\nabla}P+P\wedge P\\
& =& \pi^*R-P(\pi^*R)-P\wedge P.                    \label{eqpi}
\end{eqnarray}
Now, recall the twistor space is always a complex 2-manifold and
$D$ is a \C-linear connection. Moreover, in dimension $n=1$ we
also have that $R^\nabla$ is proportional to $\omega$ and so it is
type (1,1) for all $j$ in any fibre of the twistor space --- an
assertion equivalent to $\pi^*R$ being (1,1) for $\jnab$. On the
other hand, $P({\jnab}^+X)=\Phi^+P(X)$ so, if we prove
$[\m_{_J}^+,\m_{_J}^+]=0$, then we may conclude $R^D$ is type
(1,1). The result now follows for both vector bundles referred, by
the theorem of Koszul-Malgrange previously mentioned.

If $A,B\in\m_{_J}$, then
\[ J^+AJ^+B=AJ^-J^+B=0\]
where $J^+,J^-$ are the projections onto the $+$ or
$-\sqrt{-1}$-eigenbundles.
\end{proof}
As the reader may notice, the result is valid in any dimension
once $\jnab$ is integrable. We combine the proof above with
equation \eqref{eq4} in theorem \ref{teo1.1}.

Finally, we reach a goal: if one assigns a metric to $M$, then all
previous constructions follow and one is left with a new tool in
the study of Riemann surfaces. Letting ${\cal F}$ denote one of
the sheaves of germs of holomorphic sections of ${\cal H}^\nabla$
or ${\cal V}$, then
\[ R^1\pi_*{\cal F} \]
may tell us something new about $M$. Indeed, at the end of the
last section we discuss and conjecture that $R^1\pi_*{\cal O}$ is
non zero.

\vspace{1cm}

\section{The Penrose Transform}

Let $Z$ be a complex manifold of dimension $m$. Recall that $Z$ is
said to be strongly $q$-pseudoconvex if it admits a smooth
exhaustion function which is strongly $q$-pseudoconvex outside of
a compact subset, {\it ie.} there exists $\phi:Z\rightarrow\R$ of
class $\cinf{}$ such that the level sets $\{ x\in Z:\ \phi(x)<c
\},\ c\in\R,$ are relatively compact in $Z$, the exhaustion, and
such that the Levi form
\[  L(\phi):TZ\otimes TZ\longrightarrow\R  \]
has at least $m-q+1$ positive eigenvalues in the complement of a
compact subset $C$. If $C=\emptyset$, then $Z$ is said to be
holomorphically $q$-complete.

Recall that
\[  L(\phi)=4\sum_{i,j}\frac{\partial^2\phi}{\partial z_i\partial\overline{z}_j}
\dx z_i\otimes\dx\overline{z}_j  \] is a Hermitian 2-tensor,
independent of choice of the chart $(z_1,\ldots,z_m)$ of $Z$. From
the definition we see that $q$-completeness implies
$q+1$-completeness. Holomorphically 1-complete manifolds are known
as Stein manifolds.
\begin{teo}[\cite{Wu}]  \label{teo6.1}
A simply connected complete K\"ahler manifold $X$ of everywhere
nonpositive sectional curvature is a Stein manifold.
\end{teo}
The proof of this theorem due to H. Wu contains the following
arguments. Let $d:X\rightarrow\R$  be the Riemannian distance
function from a fixed point $p\in X$. Then it is proved that $d^2$
is smooth and strictly plurisubharmonic. It is an exhaustion
function due to completeness: a bounded and closed set is compact.

Besides $\C^n$, the canonical example to which the theorem above
applies is the Siegel domain. So the square of the distance
function in $Sp(2n,\R)/U(n)$ with invariant hyperbolic metric is
$\cinf{}$. Notice the same result does not apply to all components
of $G/U^l$, as their natural metrics may be indefinite. Yet they
are Stein spaces as we have remarked earlier.

Now let $(M,\omega,\nabla)$ be a symplectic manifold of dimension
$2n$ with a symplectic connection of Ricci type, {\it ie.} with
vanishing Weyl curvature tensor. Consider the twistor space
$(\zo,\jnab)$, which is then a complex manifold of dimension $n+k$
with $k=n(n+1)/2=\dim$ Siegel domain $G/U^0$.
\begin{lema}  \label{lema6.1}
Let $D$ be a domain in $\C^m$ and $X$ a regular complex analytic
subspace. If $\psi\in\mathrm{C}^2_D$ then
\[ L(\psi)_{|TX\otimes TX}=L(\psi_{|X}).  \]
\end{lema}
\begin{proof} We know that for every $z\in X$ there is a chart
$(z_1,\ldots,z_m)$ in a neighborhood $U$ of $z$ such that $X\cap
U=\{z\in U: z_{k+1}=\ldots=z_m=0\}$. Since $T_z(X\cap U)=\{ u\in
T_zU:\ \dx z_i(u)=0,\ i>k\}$ we find the result just by looking at
the definition of the Levi form.
\end{proof}
$M$ always admits a smooth and compatible almost complex structure
$J$, so we define a smooth function $h$ on $\zo$ to be the square
of the distance in each fibre to the section $J$ ---  which we
know to arise from a smooth Riemannian metric on the vertical
bundle $\ker\dx\pi$.
\begin{teo}  \label{teo6.2}
If $M$ has a smooth exhaustion function $\phi$, then $\zo$
is\linebreak[0] $n+1$-complete.
\end{teo}
\begin{proof} Since it is easy to prove $\phi^2$ is also an
exhaustion function, we may already assume $\phi$ to be positive.
Now let
\[ \psi= h+\phi\circ\pi  .   \]
This is a smooth and exhaustion function. To prove this notice
that $h$ is positive, so the closed level sets of $\psi$ are
inside the closed level sets of $\phi\circ\pi$ for the same
constant $c$. These
 project onto a compact subset $K_c$ of $M$. Then we have that
\begin{eqnarray*}
 \left\{ j\in\zo:\ \psi(j)\leq c\right\}& \subset
 &\left\{j\in\pi^{-1}(K_c):\ \psi(j)\leq c\right\} \\
&\subset &\bigl\{j\in\pi^{-1}(K_c):\ h(j)\leq c+\sup_{K_c}
\phi\bigr\}
\end{eqnarray*}
and, since the biggest set is compact, the closed level sets of
$\psi$ are compact.

Now, for any $x\in M$, we apply the lemma to the complex
submanifold $\pi^{-1}(x)$ and use the previous theorem to find
that $L(\psi_{|\pi^{-1}(x)})=L(h_{|\pi^{-1}(x)})$ has $k$ positive
eigenvalues. Since $L(\psi)$ is Hermitian symmetric, there is an
orthogonal complement for $\ker\dx\pi$ and we may conclude that
$L(\psi)$ has at least $k$ positive eigenvalues. Hence $\zo$ is
$q$-complete, where $q$ is such that $n+k-q+1=k$.
\end{proof}
{\ }\\
{\bf Example 1.}  If $M$ is compact we may take $\phi=0$ in the
theorem above. In particular, $\zo_{\Cpequeno\Ppequeno^n}$ is
$n+1$-complete (and not less).
\vspace{3mm}\\
{\bf Example 2.} If $M$ has some Riemannian structure for which
there is a pole, {\it ie.} there exists $x_0\in M$ such that
$\exp:T_{x_0}M\rightarrow M$ is a diffeomorphism, then we may take
$\phi=\|\exp^{-1}\|^2$.
\vspace{3mm}\\
{\bf Example 3.} Let $M=B_\epsilon(0)$, the open ball of radius
$\epsilon$ in $(\R^{2n},\omega_0)$. For this case we find
$\phi(x)=-\log(\epsilon^2-\|x\|^2)$, which is the famous function
of K. Oka.
\vspace{4mm}\\
One must realise now that the difficult thing is to find
completeness below $n+1$. This will certainly involve the
horizontal part of $\jnab$, which so much characterises twistor
spaces.
\vspace{4mm}\\
\textbf{Remark.} In general, it is impossible to find a better
result than that of the theorem: we know the Levi-Civita
connection of $\C\Pro^n$ is of Ricci type and we have seen that
parallel complex structures embed holomorphically into the twistor
space. On the other hand it is known that a $q$-complete space
does not have $n$-dimensional compact analytic submanifolds, for
any $n\geq q$.

By the same token, the {\it K\"ahlerian} twistor space ${\cal
Z}_{\Tpequeno^{2n}}$ is just holomorphically $\,n+1$-complete, and
no less.

However, with some restriction, it may well happen that it is
possible to carry on. As we have seen in example 1.1 of the
Examples section, together with theorem \ref{teo4.1}, the twistor
space of $\R^2$ with trivial connection $\nabla^0$, and hence with
all $\sigma\cdot\nabla^0$, is 1-complete or Stein (such property
is preserved under biholomorphism).

We are now ready to show the Penrose Transform.

In a parallelism with what was done in
\cite{Atiyah1,Douady,Rawnsley4}, in the celebrated Riemannian case
of $\C\Pro^3\rightarrow S^4$, we define the ``Penrose Transform''
in the symplectic context to be the direct image of any complex
analytic sheaf over twistor space onto the base manifold. Thus a
functor ${\cal O}\rightarrow\cinf{}$.
\begin{teo}  \label{teo6.3}
Let $(M,\omega,\nabla)$ be as above and ${\cal F}$ a coherent
analytic sheaf over $\zo$. Then
\[  R^q\pi_*{\cal F}=0,\ \ \ \ \ \ \ \forall q\geq n+1.  \]
\end{teo}
\begin{proof} Recall $R^q\pi_*{\cal F}$ is the sheaf associated to the presheaf
$U\mapsto H^q(\pi^{-1}U,{\cal F}).$ Hence the stalk at $x\in M$ is
\[ {\lim_{U\ni x}\mathrm{ind}}\ H^q(\pi^{-1}U,{\cal F}) .   \]
Now, for a sufficiently small neighborhood $U$ of $x$, there is a
chart $\sigma:U\rightarrow B\subset\R^{2n}$ such that
$\sigma^*\omega_0=\omega$ and $\sigma(x)=0$. Since we have a
theorem  saying there is a biholomorphism
\[ \Sigma : (\zo_U,\jnab)\longrightarrow(\zo_B,{\cal J}^{\sigma\cdot\nabla}),   \]
we may suppose our base space is $B$ and the {\it coherent}
analytic sheaf is $\Sigma_*{\cal F}$.

Finally, the $\{ B_\epsilon(0)\}_{\epsilon>0}$ form a basis for
the neighborhoods of $0$ and, by example 3 above, all
$\zo_{B_\epsilon}=\pi^{-1}(B_\epsilon)$ are n+1-complete. By
definition of inductive limit we find that $(R^q\pi_*{\cal
F})_x=0,\ \forall q\geq n+1$, appealing to Andreotti-Grauert's
``t. de finitude pour la cohomologie des espaces complexes'' ({\it
cf.} \cite{Grau}).
\end{proof}
Although we know $H^q(\pi^{-1}(x),\iota^*{\cal F})=0,\ \,\forall
q\geq 1$, where $\iota$ is the inclusion map, one has to notice in
the above proof that the $\{\pi^{-1}(U)\}$ do not form a basis of
the neighborhoods of $\pi^{-1}(x)$, as they always do in the
Riemannian case (the fibre is compact).
\vspace{4mm}\\
\textbf{Remark.} Consider example 1.1 in section 3. Recall the
global chart $(z,w)\mapsto(\xi,w)$ where $\xi=z-w\overline{z}$ is
a complex coordinate. Let us fix $z\in\R^2$ and denote
\[  X_\varepsilon=\bigl\{(\xi,w)\in\C^2:\ |w|<1,\ |z|<\varepsilon\bigr\} \]
{\it ie.} the image in $\C^2$ of $\pi^{-1}(B_\varepsilon(z))$. The
condition $|z|<\varepsilon$, where
$z=\frac{\xi+w\overline{\xi}}{1-|w|^2}$, shows that
$X_\varepsilon$ is not pseudoconvex: we can show that some regions
in the boundary are convex and others are concave. Pseudoconvexity
in $\C^n$, $n>1$, is the same as being Stein, so we may follow
\cite{Peter} and conjecture that
$H^1(X_\varepsilon,\mathcal{O})\neq 0$.

If such conjecture is true, then we can also deduce that our
sheaves $R^1\pi_*{\cal F}$ existing always over Riemann surfaces
(see section 5.2) are not in general completely trivial.

\newpage

\bibliographystyle{plain}

\end{document}